\documentclass[review,12pt,numbers,sort&compress]{elsarticle}

\usepackage{setspace}
\usepackage{subcaption}
\usepackage{multirow}
\usepackage{multicol}
\usepackage{graphicx}
\RequirePackage{amsmath}
\usepackage{bm}
\usepackage{xfrac}
\usepackage{float}
\usepackage{arydshln} 

\usepackage{amsmath,esint,amsfonts}

\usepackage[most]{tcolorbox}
\usepackage[colorinlistoftodos]{todonotes}
\usepackage{lineno}

\usepackage{lineno,hyperref}
\usepackage{cleveref} 

\journal{arXiv}

\begin{document}

\begin{frontmatter}

\title{A Non-Nested Multilevel Method for Meshless Solution of the Poisson Equation in Heat Transfer and Fluid Flow}

\author{Anand Radhakrishnan}
\author{Michael Xu}
\author{Shantanu Shahane}
\author{Surya Pratap Vanka\fnref{Corresponding Author}}
\address{Department of Mechanical Science and Engineering\\
	University of Illinois at Urbana-Champaign \\
	Urbana, Illinois 61801}
\fntext[Corresponding Author]{Corresponding Author Email: \url{spvanka@illinois.edu}}


%

\begin{abstract}
We present a non-nested multilevel algorithm for solving the Poisson equation discretized at scattered points using polyharmonic radial basis function (PHS-RBF) interpolations. We append polynomials to the radial basis functions to achieve exponential convergence of discretization errors. The interpolations are performed over local clouds of points and the Poisson equation is collocated at each of the scattered points, resulting in a sparse set of discrete equations for the unkown variables. To solve this set of equations, we have developed a non-nested multilevel algorithm utilizing multiple independently generated coarse sets of points. The restriction and prolongation operators are also constructed with the same RBF interpolations procedure. The performance of the algorithm for Dirichlet and all-Neumann boundary conditions is evaluated in three model geometries using a manufactured solution. For Dirichlet boundary conditions, rapid convergence is observed using SOR point solver as the relaxation scheme. For cases of all-Neumann boundary conditions, convergence is seen to slow down with the degree of the appended polynomial. However, when the multilevel procedure is combined with a GMRES algorithm, the convergence is seen to significantly improve. The GMRES accelerated multilevel algorithm is included in a fractional step method to solve incompressible Navier-Stokes equations.
\vspace{0.5cm}
\end{abstract}


\end{frontmatter}

\section{Introduction}

\par Multilevel methods\cite{brandt1977multi,trottenberg2000multigrid,briggs2000multigrid,yavneh2006multigrid,stuben1982multigrid} have been extensively demonstrated to provide optimal iterative convergence of discretized elliptic partial differential equations such as the Poisson equation encountered in numerous engineering problems \cite{bairstow1919two,bateman1938partial,shimada1991numerical,harlow1965numerical,poplau2001multigrid}.  By restricting residuals from fine to coarse levels, solving them on coarser levels, and prolongating corrections to the finer levels, multilevel methods provide optimal smoothing of all error frequencies encountered in an iterative solution process.  The traditional multilevel method \cite{brandt1977multi} uses nested refinement in which a finer grid is generated by subdividing a coarser grid. The refinement can be successively performed until the desired mesh fineness is obtained.  For such grids it is very easy to perform restriction and prolongation operations through distance-based interpolations.  However, if an arbitrary fine grid is first generated, then the coarsening cannot be performed by a nested procedure. In such a case, several fine grid cells (finite volumes) can be combined to make a coarser cell. Such ‘agglomeration’ multilevel procedures have been developed for Poisson equation \cite{koobus1994unstructured,chan1998agglomeration}, fpr Euler \cite{venkatakrishnan1995agglomeration,lallemand1992unstructured,bassi2012flexibility}, and Navier-Stokes \cite{langer2014agglomeration,mavriplis19963d,carre1997implicit} equations governing fluid flows. A similar technique is also implemented under the name Additive Correction Multigrid (ACM) \cite{settari1973generalization,hutchinson1988application} in which uniform corrections to the fine grid variables are obtained by solving correction equations on coarser grids generated by combining cells.  A more general technique for arbitrary sets of linear equations is the Algebraic Multigrid (AMG) method pioneered by Brandt \cite{brandt1986algebraic,ruge1987algebraic} and implemented as black-box multigrid software \cite{yang2002boomeramg,dendy1982black}.  The algebraic multigrid technique coarsens the equation system by examining the strengths of the coefficients, and hence results in efficient coarsening strategies consistent with point iterative solvers.  Lastly, non-nested multigrid techniques \cite{bittencourt2001nonnested,bittencourt2002adaptive,antonietti2019v} have also been developed as alternatives to agglomeration and additive correction multigrid methods.  The advantage of non-nested multigrid methods is that the discretization operator on coarser levels can be constructed directly by the governing equation (as in nested grids), and therefore the multigrid convergence can be optimal. However, because the grids are not nested, several independent grids need to be generated with restriction/prolongation operators constructed from multi-dimensional interpolation formulae. Recently, Katz and Jameson \cite{katz2009multicloud} used multiquadrics to discretize the compressible flow equations and accelerated the convergence to steady state by a multilevel meshless method. The coarse sets in their method were generated by a special coarsening algorithm.

\par The principal characteristic of a multigrid technique is grid-independent rate of convergence of an iterative solver for any arbitrary size of the discrete problem. The number of iterations to reach a specified level of convergence depends only on the rate of convergence of the low frequency errors on the coarsest level. This is a significant advantage of multigrid based numerical algorithms for elliptic problems. Multigrid methods can also be considered as preconditioners to iterative solvers \cite{wienands2000fourier,evstigneev2019numerical,elman2001multigrid,erlangga2006novel,xu1996auxiliary}, lowering the condition number of the iteration matrix to achieve fast convergence. Multigrid methods have been used for solving sets of coupled equations such as those of fluid flow employing either decoupled or coupled solution of the Navier-Stokes equations \cite{vanka1986block,john2000numerical,thompson1989adaptive,paisley1999multigrid,sivaloganathan1988multigrid}. For external flows on unstructured grids, agglomeration multigrid with Runge-Kutta time marching schemes have been developed \cite{mavriplis1990multigrid}.

\par A recent trend in discretizing partial differential equations in complex domains has been the concept of meshless techniques that use scattered points to represent a complex domain \cite{hardy1971multiquadric,kansa1990multiquadrics,kansa1990multiquadrics2,franke1998solving,fasshauer1996solving,shu2003local,ding2007numerical,sanyasiraju2008local,wong2002compactly,shahane2020high}.  In such meshless methods, a variable is first interpolated either globally \cite{kansa1990multiquadrics} or locally over a “cloud” \cite{kansa1990multiquadrics2} of nearby scattered points using weighted sums of a basis function. The variable is subsequently collocated at the scattered points to evaluate the interpolation coefficients. If the values of any variable are known at the scattered points, the interpolation can then be used to calculate values at in-between locations within the cloud. However, if the values are unknown but satisfy an underlying differential equation, the governing equation can be satisfied discretely at the scattered points either by collocation or by the method of weighted residuals \cite{franke1998convergence,yun2016improved,bourantas2019explicit,yang2019solving,varanasi2010meshless}. The resulting set of linear (or nonlinear) equations is then solved for the unknown values of the dependent variable at the scattered points.

\par Some of the common RBFs previously considered for solution of partial differential equations are as follows:
\begin{equation}
	\begin{aligned}
		\text{Multiquadrics (MQ): } \phi(\bm{r})=&(\bm{r}^2 + c^2)^{1/2}\\
		\text{Inverse Multiquadrics (IMQ): } \phi(\bm{r})=&(\bm{r}^2 + c ^2)^{-1/2}\\
		\text{Gaussian: } \phi(\bm{r})=&\exp\left(\frac{-\bm{r}^2}{\sigma ^2}\right)\\
		\text{Polyharmonic Splines (PHS): } \phi(\bm{r})=&\bm{r}^{2a+1},\hspace{0.1cm} a \in \mathbb{N}\\
		\text{Thin Plate Splines (TPS): } \phi(\bm{r})=&\bm{r}^{2a} log(\bm{r}),\hspace{0.1cm} a \in \mathbb{N}\\
	\end{aligned}
	\label{Eq:RBF_list}
\end{equation}
where $r$ is the magnitude of the distance between two points. The first three radial basis functions need specification of a shape factor c or $\sigma$ which controls the shape of the RBF. A variable $f(\bm{x})$ is interpolated as:

\begin{equation}
	f(\bm{x}) = \sum_{i=1}^{N} \alpha_i \phi(||\bm{x} - \bm{x_i}||_2)
	\label{Eq:RBF_interp}
\end{equation}

where the symbol $N$ denotes the number of points used in the interpolation ‘cloud’ for point $\bm{x}$
and $\alpha_{i}$
denote the weighting coefficients. In global interpolation methods, the value of $N$ is equal
to the total number of discrete points in the domain. Global methods, while highly accurate,
produce interpolation matrices that can have large condition numbers even for modest values of
$N$. Therefore, they are not convenient to solve problems with large numbers of scattered points.
The shape parameter influences the accuracy of the interpolation and the condition number of
the matrix to evaluate the interpolation coefficients. The prescription of the optimum value of
these shape parameters has been one of the difficult issues in the robust use of RBFs for
interpolation \cite{buhmann2003radial,franke1980smooth}.

\par Interpolation within a finite number of points in a cloud produces a sparse assembled matrix with
a manageable condition number. However, the condition number can still become large for large
cloud sizes and large total number of points. Further, with multiquadric, inverse multiquadric and
Gaussian interpolations \cite{fornberg2011stable} a ‘saturation problem’ can occur, which prevents the discretization
errors to asymptotically converge to machine precision. Interpolation with polyharmonic splines does not need the additional shape parameter but also encounters the saturation
problem. However, if the PHS-RBF is appended with polynomials \cite{flyer2016role,shahane2020high}, the discretization error can be
made to decrease exponentially as the degree of the appended polynomial. Also,
interpolation using PHS with an appended polynomial converges the discretization errors to
machine precision, although appending a polynomial increases the cloud size and the
computational work. The coefficients of the PHS-RBF interpolation and of the appended
monomials are evaluated by collocating the variable at the scattered locations supplemented
with constraints on the appended monomials. The high order of accuracy and convergence to
machine precision makes the PHS-RBF interpolation with an appended polynomial an attractive
meshless method \cite{bayona2017role}. \nocite{barnett2015robust,miotti2021fully}

\par In this paper, we describe a multilevel procedure to efficiently solve the Poisson equation in
complex domains discretized by the PHS-RBF method. The Poisson equation is encountered while
solving numerous engineering problems such as heat conduction, incompressible flows,
electromagnetics, porous media flows, etc. Many problems involve irregular solution domains
whose accurate discretization and efficient solution are of much importance. Complex domains
have been traditionally discretized using curvilinear \cite{anderson1995computational,pletcher2012computational} or unstructured \cite{fletcher2012computational,mathur1997pressure} finite difference or finite volume methods. For spectral accuracy, the
finite element method has also been combined with Chebyshev polynomial expansions, giving
rise to the well-known spectral element method \cite{fischer1989parallel,xu2018spectral,karniadakis2013spectral}. The spectral element method has high
accuracy but requires the placement of the points within an element at the specified roots of
the polynomials. While the discretization convergence of the finite difference method depends
on the order of truncation of the Taylor series, the accuracy of the finite volume method depends
on the evaluation of the interface fluxes.

\par The PHS-RBF discretization can also be formulated to control the local accuracy by selectively
varying the degree of the appended polynomial and the spacing between points. This is akin to
($h-p$) refinement used in higher order finite element and spectral element methods.
Discretization of the Poisson equation using the PHS-RBF interpolants over finite sized clouds of
points results in a sparse set of linear equations with sparsity proportional to the highest degree
of the appended monomial. The condition number of the coefficient matrix usually becomes
large, requiring powerful iterative solvers for fast convergence. Conventional single level solvers (such as Jacobi, and SOR) have slow convergence for the long
wavelength errors and converge slowly. This motivates us to investigate convergence of multilevel methods as accelerators with
low-cost solvers .

\par The meshless method does not have any underlying grid, and hence is not amenable to
traditional nested coarsening or agglomeration of cells. Black-box multigrid procedures based on
AMG use heuristic techniques to generate the coarse set of equations. However, this approach has not yet been explored for meshless methods. Non-nested multilevel procedures in which
multiple coarse levels of points are generated either independently or by grouping a few surrounding
points into a “mean” point show promise to interpolate variables, residuals, and corrections
across levels. To our knowledge, such multilevel methods have not been explored for the PHS-
RBF solution of elliptic equations in complex domains. For such non-nested multigrid methods,
the interpolation across the various coarse and fine levels of points can be performed using the
same radial basis functions also used on the finest set of points. In this work, we describe such a
meshless multilevel method for solving the Poisson equation and evaluate its convergence in
number of model problems with Dirichlet and all-Neumann boundary conditions. The
performance of the method is evaluated for an oscillatory manufactured solution of varying wave
number. For all-Neumann boundary conditions the equations are regularized by adding a
constraint that makes the linear system well-conditioned \cite{regularize2019}.

\section{The PHS--RBF Method} \label{Sec:PHS--RBF Method}
The PHS-RBF interpolates a function $f(\bm{x})$ whose values are known at scattered points using weighted kernels as:

\begin{equation}
	f(\bm{x}) = \sum_{i=1}^{N} \alpha_i \phi(||\bm{x} - \bm{x_i}||_2)
	\label{Eq:RBF_interp}
\end{equation}

where $\phi$ is the polyharmonic spline function given as:

\begin{equation}
	\phi(||\bm{x} - \bm{x_i}||_2)  = ||\bm{x} - \bm{x_i}||_2^{2p+1}
	\label{Eq:PHS}
\end{equation}

and $p$ = 1,2,3, etc. Further, as mentioned earlier, PHS are appended with number of monomials,
the number given by $\binom{l}{d}$,
where $l$ is the desired highest degree of the monomial and $d$ is the dimension of the problem. Thus, \cref{Eq:RBF_interp} is extended as:

\begin{equation}
	f(\bm{x}) = \sum_{i=1}^{N} \alpha_i \phi(||\bm{x} - \bm{x_i}||_2) + \sum_{j=1}^{M} \gamma_j P_j(\bm{x})
	\label{Eq:RBF_PHS_interp}
\end{equation}

The RBF $\phi(||\bm{x} - \bm{x_i}||_2)$ is a scalar function of the Euclidian distance between the points
irrespective of the problem dimension. $P_j(\bm{x})$ denotes the $j^{th}$ degree monomial which is weighted
by the coefficient $\gamma_{j}$.
The $(N+M)$ unknown coefficients are evaluated by collocating the variables
at the discrete points and satisfying the constraints on the coefficients of polynomials.

\begin{equation}
    \sum_{i=1}^{N} \alpha_i P_j(\bm{x_{i}}) = 0; \hspace{0.5cm}   j = 1,2,...,M
	\label{Eq:PHS_Colloc}
\end{equation}

Usually, the total number of cloud points $N$
is selected to be a multiple (around twice) of the
number of monomials. The collocation equations and constraints for $\bm{\alpha}$ and $\bm{\gamma}$ can be written as:

\begin{equation}
	\begin{bmatrix}
		\bm{\Phi} & \bm{P}  \\
		\bm{P}^T & \bm{0} \\
	\end{bmatrix}
	\begin{bmatrix}
		\bm{\alpha}  \\
		\bm{\gamma} \\
	\end{bmatrix} =
	\begin{bmatrix}
		\bm{f}  \\
		\bm{0} \\
	\end{bmatrix}
	\label{Eq:RBF_interp_mat_vec}
\end{equation}

or,

\begin{equation}
	\begin{bmatrix}
        \bm{D} \\
	\end{bmatrix}
	\begin{bmatrix}
		\bm{\alpha}  \\
		\bm{\gamma} \\
	\end{bmatrix} =
	\begin{bmatrix}
		\bm{f}  \\
		\bm{0} \\
	\end{bmatrix}
	\label{Eq:RBF_interp_mat_vec2}
\end{equation}

where the superscript $T$ denotes the transpose, $\bm{\alpha} = [\alpha_1,...,\alpha_N]^T$, $\bm{\gamma} = [\gamma_1,...,\gamma_M]^T$, $\bm{f} = [f(\bm{x_1}),...,f(\bm{x_N})]^T$ and $\bm{0}$ is the matrix of all zeros of appropriate size. Sizes of the submatrices $\bm{\Phi}$ and $\bm{P}$ are $N\times N$ and $N\times M$ respectively. The matrix $\Phi$ is given by,

\begin{equation}
	\bm{\Phi} =
	\begin{bmatrix}
		\phi \left(||\bm{x_1} - \bm{x_1}||_2\right) & \dots  & \phi \left(||\bm{x_1} - \bm{x_q}||_2\right) \\
		\vdots & \ddots & \vdots \\
		\phi \left(||\bm{x_q} - \bm{x_1}||_2\right) & \dots  & \phi \left(||\bm{x_q} - \bm{x_q}||_2\right) \\
	\end{bmatrix}
	\label{Eq:RBF_interp_phi}
\end{equation}

For a polynomial of maximum degree two, the matrix is given as:

\begin{equation}
	\bm{P} =
	\begin{bmatrix}
		1 & x_1  & y_1 & x_1^2 & x_1 y_1 & y_1^2 \\
		\vdots & \vdots & \vdots & \vdots & \vdots & \vdots \\
		1 & x_q  & y_q & x_q^2 & x_q y_q & y_q^2 \\
	\end{bmatrix}
	\label{Eq:RBF_interp_poly}
\end{equation}

Given the values of $\bm{f}$ at the scattered points, the values of $\bm{\alpha}$ and $\bm{\gamma}$ can be evaluated by solving \cref{Eq:RBF_interp_mat_vec2}. The function $\bm{f}$
can then be evaluated at any arbitrary location $\bm{x}$
within the cloud. These
relations can be implicitly used to determine the discrete values of $\bm{f}$ that satisfy a given governing
equation at the discrete locations.

\par To determine the unknown values of the variable, the governing differential equation is
collocated at the scattered locations \cite{kansa1990multiquadrics,kansa1990multiquadrics2}. The derivatives at any location are evaluated
by differentiating the analytical interpolation functions and satisfying the governing partial differential equation at every discrete point, thus giving the desired number of equations to
determine the values of $\bm{f}$. Let $\mathcal{L}$ denote the Laplacian ($\nabla^{2}$) operator. Then, we can write:

\begin{equation}
	\mathcal{L}[{f(\bm{x})}] =
	\begin{bmatrix}
		\mathcal{L}[\bm{\Phi}] & \mathcal{L}[\bm{P}]  \\
	\end{bmatrix}
	\begin{bmatrix}
		\bm{\alpha}  \\
		\bm{\gamma} \\
	\end{bmatrix}
	\label{Eq:RBF_interp_mat_vec_L}
\end{equation}

\begin{equation}
    =
	\begin{bmatrix}
		\mathcal{L}[\bm{\Phi}] & \mathcal{L}[\bm{P}]  \\
	\end{bmatrix}
	\begin{bmatrix}
        \bm{D} \\
	\end{bmatrix}^{-1}
	\begin{bmatrix}
		\bm{f} \\
		\bm{0} \\
	\end{bmatrix}
	=
	\begin{bmatrix}
        \bm{B}_{1} & \bm{B}_{2} \\
	\end{bmatrix}
	\begin{bmatrix}
		\bm{f} \\
		\bm{0} \\
	\end{bmatrix}
	=
	\begin{bmatrix}
        \bm{B}_{1} \\
	\end{bmatrix}
	\begin{bmatrix}
        \bm{f} \\
	\end{bmatrix}
	\label{Eq:RBF_interp_mat_vec_L_solve}
\end{equation}

The discrete coefficients for estimation of the Laplacian operator at each point are then
combined in a global matrix of $n$ rows (where $n$ is the total number of scattered points), resulting
in a set of linear equations given by,

\begin{equation}
	\begin{bmatrix}
        \bm{A} \\
	\end{bmatrix}
	\begin{bmatrix}
        \bm{X} \\
	\end{bmatrix}
	=
	\begin{bmatrix}
        \bm{b} \\
	\end{bmatrix}
    \label{Eq:Sparse_eq}
\end{equation}

where $\bm{A}$ is a sparse matrix of the coefficients of size equal to the number of the scattered points and $\bm{b}$ is the vector of
the source terms. $\bm{X}$ is the vector of all the unknown discrete values. The matrix $\bm{A}$ depends only on the coordinates of the cloud points and can be computed and stored at the beginning of
the computational algorithm. The computation of the $\bm{A}$ matrix requires solution of dense linear
systems corresponding to the cloud of each discrete point. However, since the computations for an
individual point are independent of others, the coefficients can be computed in parallel. Further
economies can be gained by considering groups of points as one cluster and using the same $\bm{D}$
matrix for points in the cluster \cite{shankar2017overlapped}.

\par Any traditional iterative scheme can be used to solve \cref{Eq:Sparse_eq}. However, the coefficients in $\bm{A}$ are highly oscillatory, making traditional solvers slow as well as difficult to converge. Further, single grid iterative solvers are inefficient when the matrix size is large. We
describe below a multilevel procedure utilizing non-nested coarsening and appropriate
restriction and prolongation procedures.
\section{Details of the Multilevel Solution Procedure}

Multilevel iterative procedures primarily consist of four main components. First, it is necessary
to define coarse levels for solution of the equations for corrections. In hierarchical refinement
of structured grids, the fine and coarse grids are nested, and the fine grids are defined by
sequentially halving the inter-point distances or subdividing the elements in case of finite
element or finite volume methods. Alternatively, the coarse grids can also be defined by
progressively agglomerating the finite volumes defined on the finest grid. We use here the concept of non-nested multilevel
methods where the individual grids of desired fineness are generated independently. In the
context of meshless methods, we use multiple non-nested sets consisting of different numbers
of scattered points. The second component of multilevel schemes is the relaxation operator (iterative scheme) used
to solve the appropriate discrete equations on any given level. Such a relaxation scheme can
range from a point solver (Jacobi, Gauss-Seidel, SOR), to more powerful Krylov sub-space based
solvers such as conjugate gradient, BiCGSTAB and GMRES \cite{saad1986gmres}. The third component of multilevel methods is the interpolation procedure to transfer residuals
(restriction) from finer levels to coarser level and corrections (prolongation) from coarser level to finer level. In geometry-based multilevel
methods, restriction and prolongation operators are constructed by distance weighting using
bilinear or trilinear interpolation.
In non-nested multigrid methods, the restriction and
prolongation require interpolation across scattered points. Finally, a cycling strategy that defines the sequence in which the fine and coarse
levels are visited needs to be specified. On any level the number of relaxations can
be fixed or adaptively controlled based on the convergence rate. Details of these four steps in our algorithm are given below.

\subsection{Generation of Coarse Level Point Sets}

As mentioned earlier, the meshless discretization only uses scattered points with no explicit
connectivity. The points can be generated by many procedures of choice, including by
conventional finite element mesh generation techniques. A straight-forward way to generate
scattered points is to use an established finite element mesh generator and extract only the coordinates of the vertices. In meshless techniques, because of the absence of connectivity
between the points, quality of the grid as characterized by the slenderness of the elements does
not explicitly affect accuracy of discretization. However, the inter-point spacings and the total
number of points influence the condition number and solvability of the system of equations. In
the present study, we generated the point sets on multiple levels using an open-source finite
element grid generator, Gmsh \cite{geuzaine2009gmsh}.

\subsection{Relaxation Operator}
The coefficient matrix $\bm{A}$ in \cref{Eq:Sparse_eq} is computed and stored in sparse row column (CSR) format
consistent with the Eigen library \cite{eigenweb}. At Dirichlet boundaries, where the values are known, the
interior equations are condensed with substitution of the boundary values. For Neumann
boundaries, special clouds for the boundary points are defined which contain only interior points
and no points on the boundary. The discrete equations of the boundary points are then
substituted into the equations of the interior points, resulting in an implicit coupling of the
boundary and interior points. The matrix of coefficients is then re-ordered using the Reverse
Cuthill-McGee (RCM) \cite{cuthill1969reducing} algorithm to reduce bandwidth and increase computational efficiency.
Currently we use the SOR point solver with an over-relaxation factor of 1.4 as the relaxation
scheme. The SOR is the most efficient of the point solvers, although the Jacobi solver is better
suited for parallelization. Other solvers such as ILU with drop tolerance, preconditioned
BiCGSTAB and GMRES can also be used along with the multilevel method. In Section 4.4 we
evaluate a combination of GMRES and the multilevel algorithm.

\subsection{Restriction and Prolongation Operators}
For the linear Poisson equation, a correction multilevel scheme is appropriate. Only the
residuals need to be restricted from a fine level to a coarser level. The residuals on the finer level
are calculated after a prescribed number of relaxations as:
\begin{equation}
    \bm{R^{h}} = \bm{b^{h}} - \bm{A}^{h}\bm{X^{h}}
    \label{Eq:Residual}
\end{equation}

where the superscript $h$ refers to the set of points on the finer level and $\bm{R}$, $\bm{b}$, and $\bm{X}$
are vectors that
denote the residual, the right-hand side, and the solution respectively. The residuals at the
coarser set of points, denoted with superscript $H$
are obtained by applying the interpolation operator $I_{h}^{H}$
which restricts the residuals as:

\begin{equation}
    \bm{R}^{H} = \bm{I}_{h}^{H} \bm{R}^{h}
    \label{Eq:Restriction}
\end{equation}

\begin{figure}[H]
	\centering
	\begin{subfigure}[t]{0.48\textwidth}
		\includegraphics[width=\textwidth]{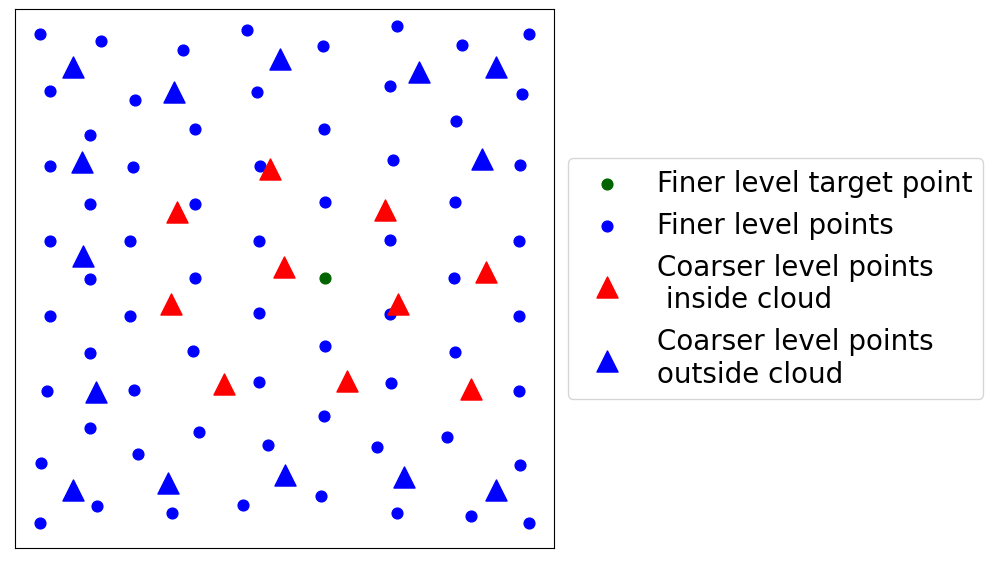}
		\caption{Restriction from finer level to coarser level}
	\end{subfigure}
	\begin{subfigure}[t]{0.48\textwidth}
		\includegraphics[width=\textwidth]{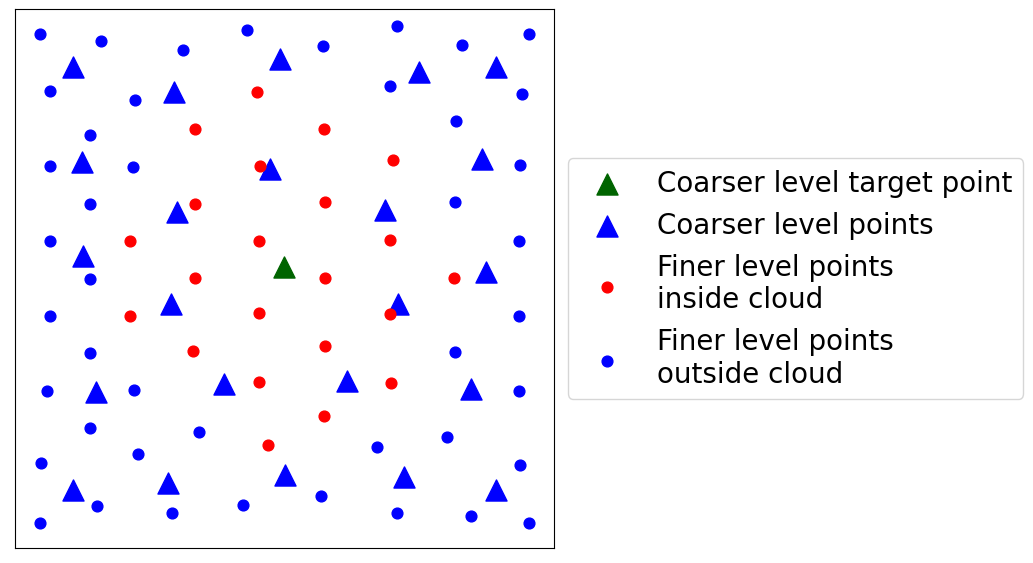}
		\caption{Prolongation from coarser level to finer level}
	\end{subfigure}
	\caption{PHS-RBF interpolation between finer and coarser levels on a square geometry}
	\label{Fig:Restrict_Prolong}
\end{figure}

The interpolation operator $\bm{I}_{h}^{H}$
is constructed as follows. Given the coordinates of the coarse and
fine sets of points, we first identify a cloud of closest fine level points for each one of the coarse
set points. The number of fine level points for each coarse set point depends on the desired
accuracy of the interpolation. We use the same PHS-RBF interpolation procedure with appended
polynomial as described earlier and  a graphical representation of this
interpolation operation is given in \cref{Fig:Restrict_Prolong}.

The order of accuracy for restricting residuals can be the same
as the order of accuracy used for discretization of the partial derivatives in the governing
equation or it can be lower. To construct the restriction interpolation matrices, we start with the
interpolation of fine level points over the selected cloud, as:

\begin{equation}
	f(\bm{x^{h}}) = \sum_{i=1}^{N^{h}} \alpha_i \phi(||\bm{x^{h}} - \bm{x_i^{h}}||_2) + \sum_{j=1}^{M^{h}} \gamma_j P_j(\bm{x^{h}})
	\label{Eq:RBF_PHS_interp_grid}
\end{equation}

where $N^{h}$ spans the cloud of fine points for a given coarse level point, and $M^{h}$ represents the number of appended monomials. Collocation of the fine set values at the discrete points of interpolation gives the coefficients $\bm{\alpha}$ and $\bm{\gamma}$. The function value is then evaluated at the
coordinates of the coarse point. The resulting weights multiplying the fine set values to the
coarse set points are evaluated at the beginning of the solution procedure and stored in a sparse
rectangular matrix. Such matrices are computed for every pair of fine and coarse levels. We
observed that the order of accuracy of interpolations need not be as high as the discretization
accuracy used on the finest level as it did not noticeably affect the rate of the multilevel
convergence. Using a lower order polynomial in the PHS-RBF interpolation matrix however
saves some CPU time. At the boundary points, the residuals are implicitly zero because the boundary equations are substituted in the interior equations. The boundary equations are therefore not solved explicitly. After the boundary equations are
eliminated, the equations for the interior points are reordered using the RCM ordering.

\par When the down-leg is completed, the corrections are successively prolongated to the points on
the finer levels. For prolongation, we again compute the interpolation matrices by first defining
the cloud of coarse set points to be used in the interpolation for each point on a finer level. The
prolongation matrices are computed in the same way as the restriction matrices by first
determining the PHS-RBF coefficients in the coarse point interpolation and then determining the
weights of each coarse point for a given fine level point. Thus, each fine level value is corrected
as:
\begin{equation}
    \bm{X}^{h} =\bm{X}^{h} + \bm{I}_{H}^{h} \delta \bm{X}^{H}
    \label{Eq:Restriction}
\end{equation}

For Dirichlet boundary conditions, the corrections on the boundaries are zero. However, for
Neumann conditions, the corrections at the boundary points are evaluated from the interior
coarse grid corrections using the discretized Neumann condition which was earlier
substituted into the interior equations. Prolongations are also followed by a relaxation step.
Both the restriction and prolongation matrices are calculated upfront of the iterations. Note
that the current restriction and prolongation operators are based solely on the geometric
coordinates of the points and do not take into accounts the strengths of the coefficients, as in
AMG.

\section{Results}
\subsection{Problems Considered}
The performance of the above algorithm is evaluated for the solution of the Poisson equation in
three model geometries using a manufactured solution. A manufactured solution satisfies a
governing equation with appropriate additional source terms obtained by substituting the
manufactured solution in the governing equation. We consider the manufactured solution:

\begin{equation}
    T(x,y) = \cos(k \pi x) \cos(k \pi y)
    \label{Eq:Manufac_Soln}
\end{equation}

which satisfies the Poisson equation,
\begin{equation}
    \nabla^{2} T(x,y) = -2\pi^{2} k^{2}\cos(k \pi x) \cos(k \pi y)
    \label{Eq:Poisson}
\end{equation}

where $k$ is the wavenumber of the manufactured solution. To test the robustness of the
multilevel convergence, we varied a) geometry; b) wavenumber; c) type of boundary condition
(Dirichlet and Neumann); d) number of scattered points on the finest level; e) number of levels
in the multilevel cycling and f) the degree of the highest appended monomial in the PHS-RBF
discretization. For Dirichlet conditions, the exact analytical values are prescribed on the
boundaries, while the exact analytical derivatives are prescribed for the Neumann boundary
condition. The effect of the order of the polynomial for constructing coarse level restriction /
prolongation and discretization operators on the multilevel performance was also investigated.
We found that using a lower order polynomial for interpolation and discretization on the coarser
levels did not noticeably affect the multilevel performance but provided a slightly reduced
computational work. The use of lower order polynomial on the coarse levels does not influence
the accuracy of the final solution obtained on the finest level.

\begin{figure}[H]
	\centering
	\begin{subfigure}[t]{0.32\textwidth}
		\includegraphics[width=\textwidth]{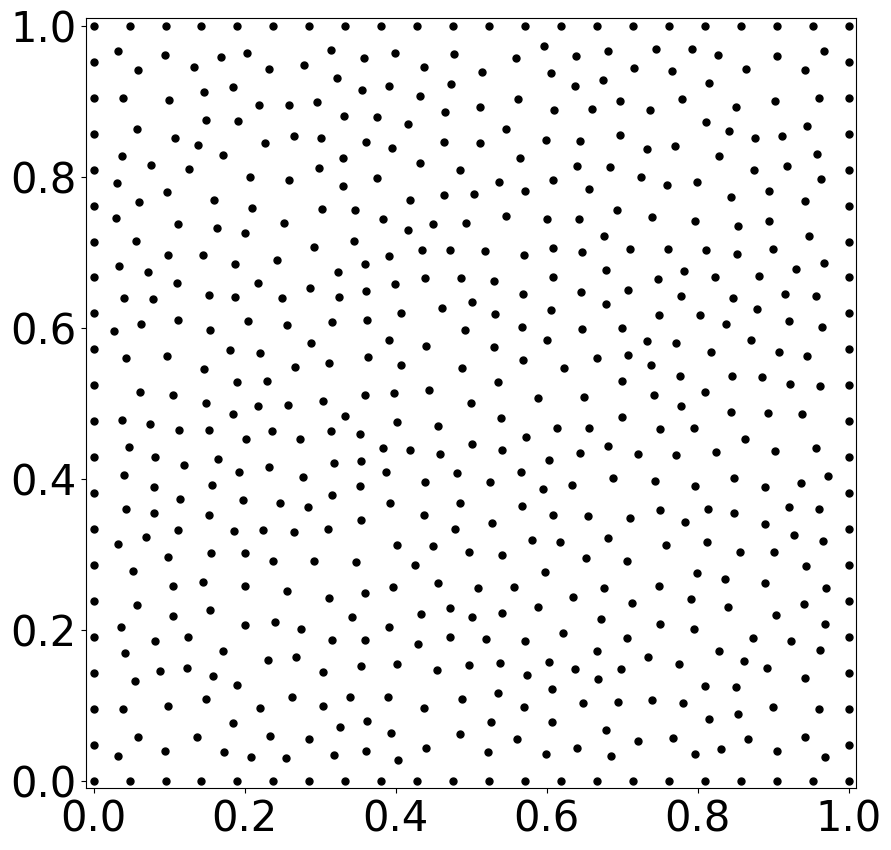}
		\caption{Square}
	\end{subfigure}
	\begin{subfigure}[t]{0.32\textwidth}
		\includegraphics[width=\textwidth]{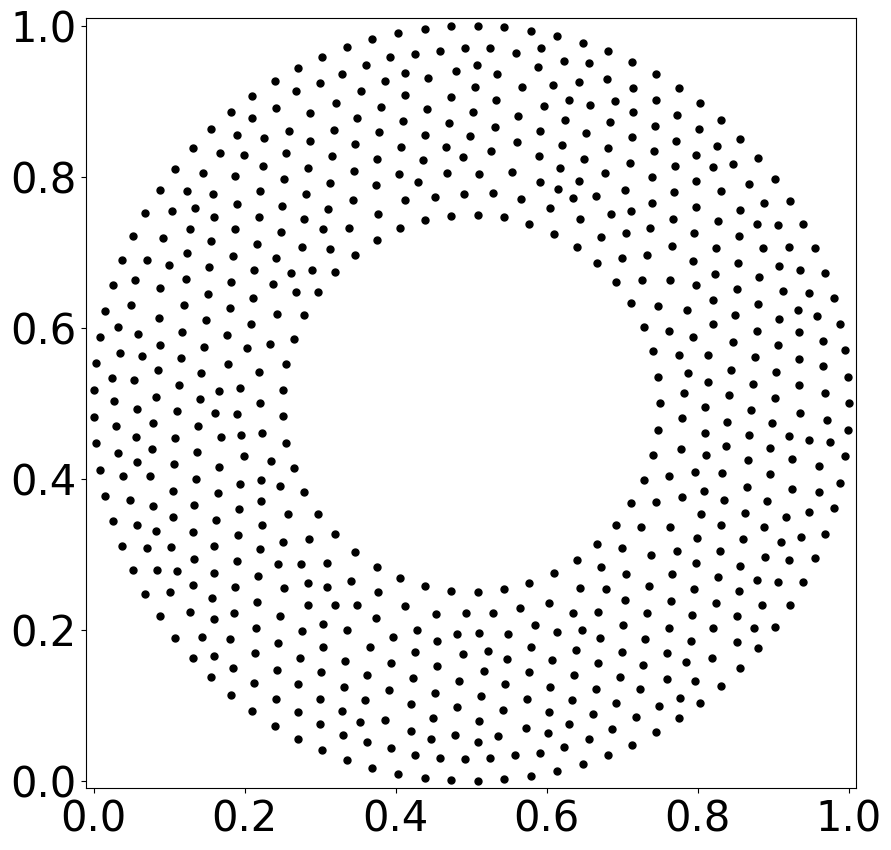}
		\caption{Circular annulus}
	\end{subfigure}
	\begin{subfigure}[t]{0.32\textwidth}
		\includegraphics[width=\textwidth]{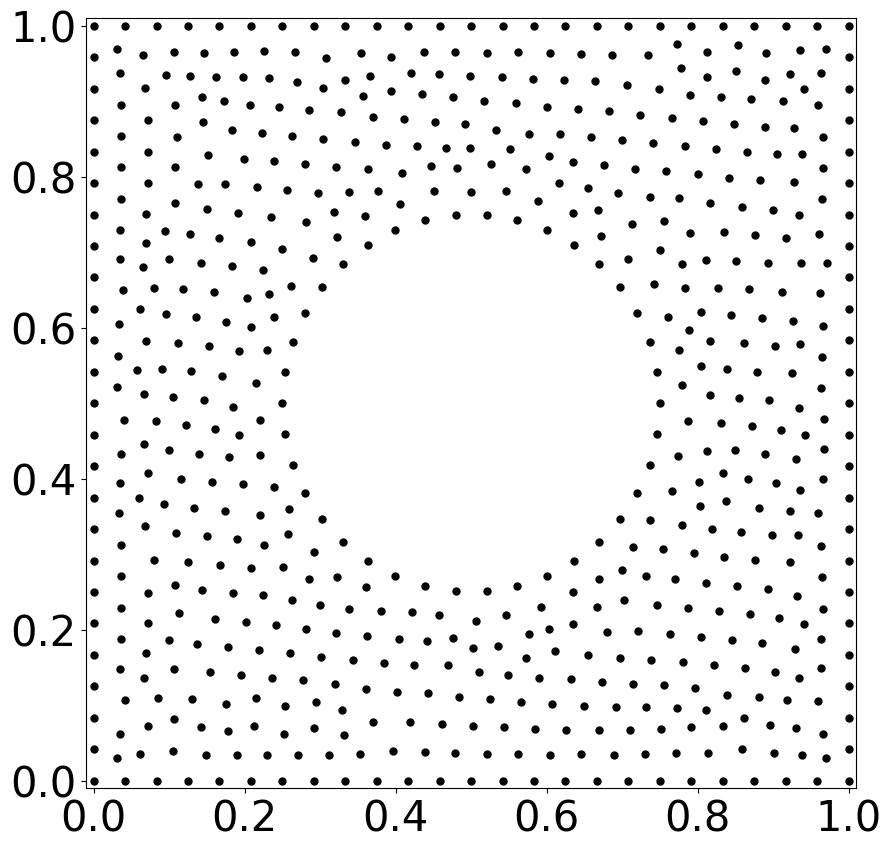}
		\caption{Square-with-hole}
	\end{subfigure}
	\caption{Geometries considered and layout of discrete points}
	\label{Fig:Scatter_Plots}
\end{figure}

\Cref{Fig:Scatter_Plots} shows the three geometries considered in this study. The scattered points shown were
obtained from vertices of triangles generated by an open-source mesh generation software \cite{geuzaine2009gmsh}. For each geometry, sets with different numbers of scattered points were obtained using the same software.

\begin{table}[H]
	\centering
	\begin{tabular}{|c|c|c|c|}
		\hline
		Geometry/Level & Square & Concentric  annulus & Square-with-hole \\ \hline
		1 & 98 & 90 & 89 \\ \hline
		2 & 169 & 188 & 176 \\ \hline
		3 & 607 & 650 & 640 \\ \hline
		4 & 2535 & 2581 & 2532 \\ \hline
		5 & 10023 & 10207 & 10197 \\ \hline

	\end{tabular}
	\caption{Numbers of points for each problem on different levels}
	\label{tab:Levels}
\end{table}

\Cref{tab:Levels} gives the numbers of points used for each problem on different levels. The discretization accuracy was computed as the $L_{1}$ norm of the difference between the converged discrete solution and the
analytical solution defined as:

\begin{equation}
    E = \frac{\sum_{i=1}^{n}|T_{PHS} - T_{exact}|_{i}}{\sum_{i=1}^{n}|T_{exact}|_{i}}
    \label{Eq:Accuracy}
\end{equation}

Convergence of
the equations was monitored by the $L_{1}$
norm of the residuals defined as:

\begin{equation}
    E = \frac{\sum_{i=1}^{n}|\nabla^{2}T_{PHS} - RHS|_{i}}{\sum_{i=1}^{n}|RHS|_{i}}
    \label{Eq:Accuracy}
\end{equation}

The multilevel iterations were continued until the residual norm reached a value below $10^{-10}$
. The initial values for the variables were always set to zero, and the residuals were normalized with the
norm of the source terms. At any iteration, the residuals were evaluated at the top of the V-
cycle. At each level, the residuals were relaxed with five SOR iterations after which they were
restricted to the next coarser level. Five SOR relaxations were also performed after prolongating the
corrections to a finer level. An over-relaxation parameter of 1.4 was used on all levels.

\subsection{Results for Dirichlet Boundary Conditions}

\begin{figure}[H]
	\centering
	\begin{subfigure}[t]{0.32\textwidth}
		\includegraphics[width=\textwidth]{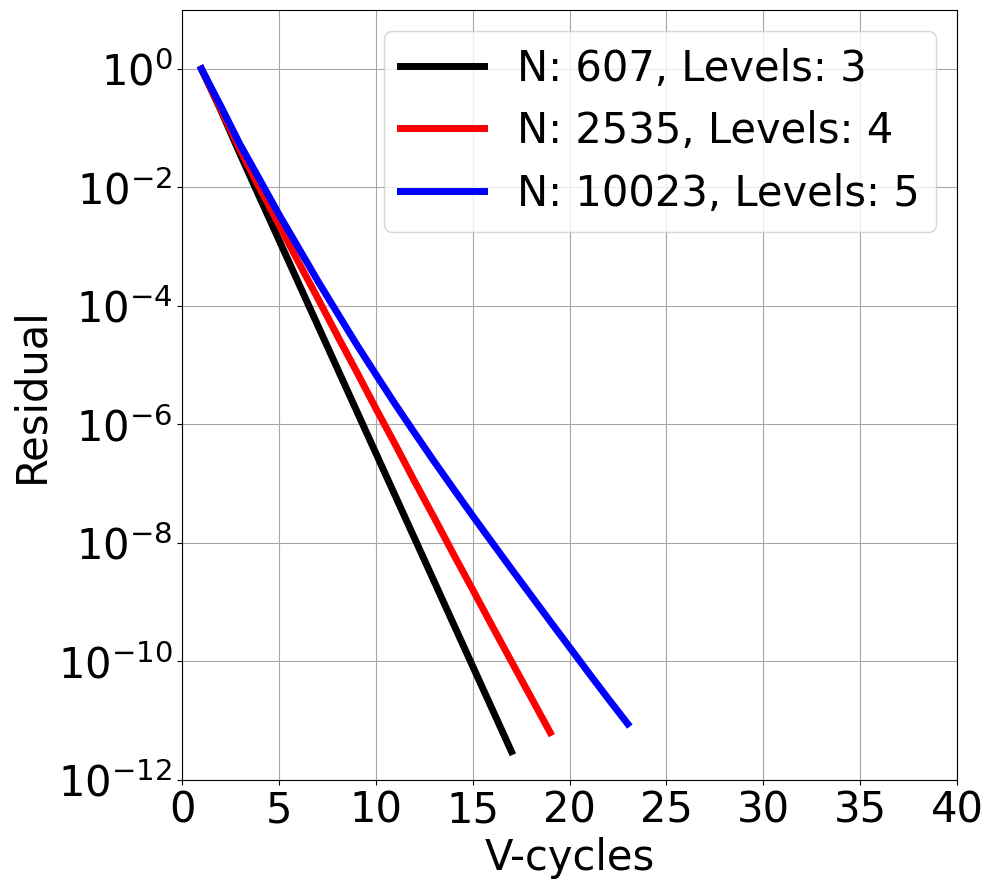}
		\caption{Degree of appended polynomial = 3}
	\end{subfigure}
	\begin{subfigure}[t]{0.32\textwidth}
		\includegraphics[width=\textwidth]{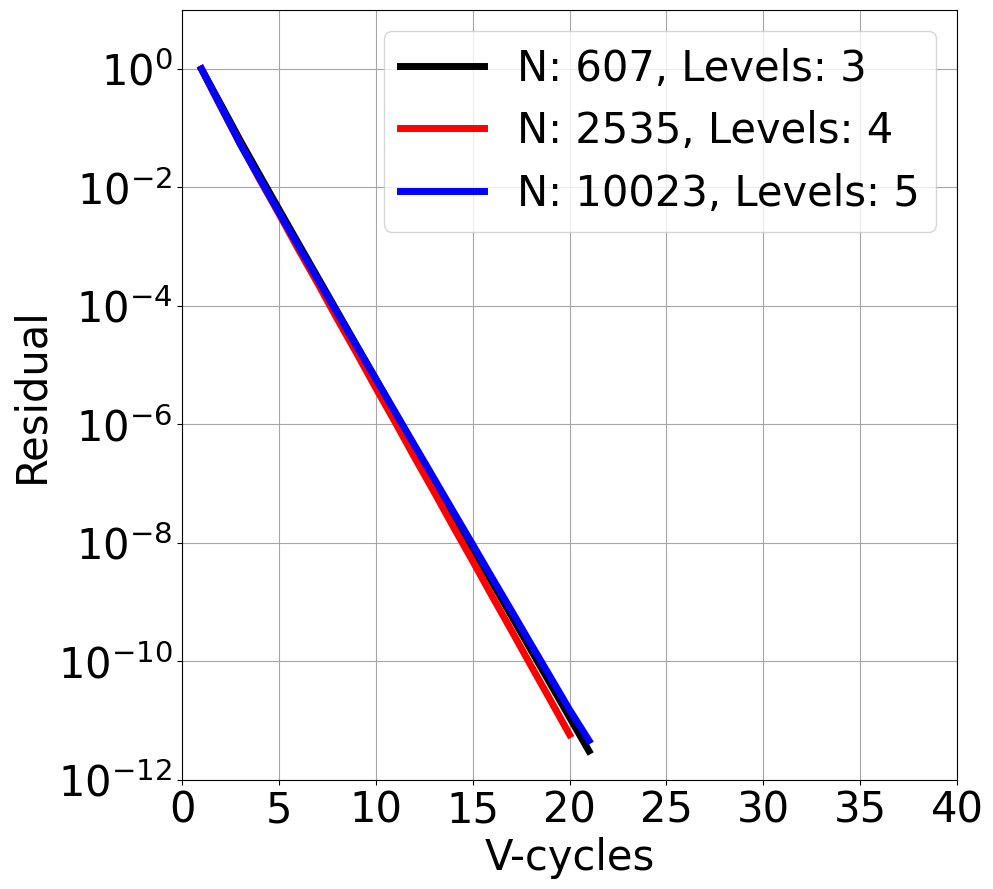}
		\caption{Degree of appended polynomial = 4}
	\end{subfigure}
	\begin{subfigure}[t]{0.32\textwidth}
		\includegraphics[width=\textwidth]{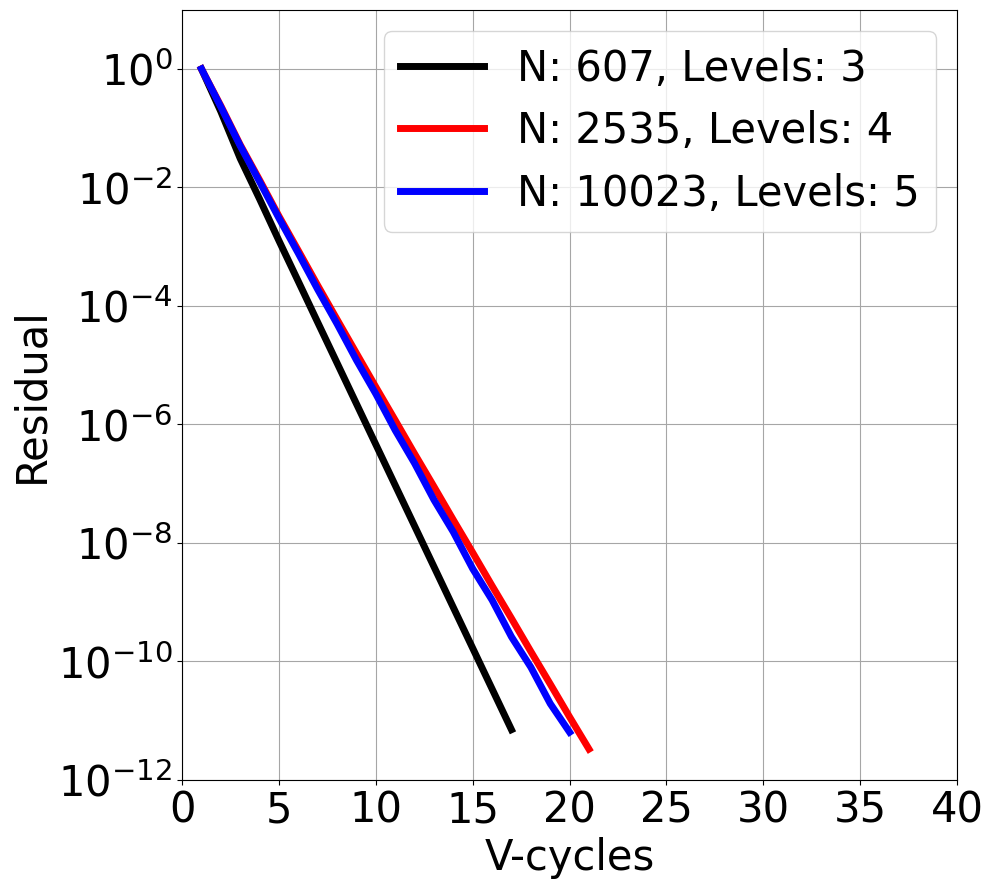}
		\caption{Degree of appended polynomial = 6}
	\end{subfigure}
	\caption{Convergence of the residual on a square with Dirichlet boundary conditions}
	\label{Fig:Dirichlet_Square}
\end{figure}

Systematic calculations were performed for each geometry by varying the highest degree of the
appended monomials, number of scattered points on the finest level, and the wavenumber of the manufactured solution. The highest degree of the appended monomials was varied between
3 and 6, and five levels of points were considered with the largest set consisting around 10,000
points. The wavenumber was varied between 1 and 4. \Cref{Fig:Dirichlet_Square} shows selected multilevel
convergence histories for the square geometry using wavenumber $k$ = 1 and three sets of points with
polynomial degree of 3, 4, and 6. For each calculation, the coarsest level had 98 points. The
number of points on all levels are given in \cref{tab:Levels}. It is seen that the residual is quickly reduced
by about ten orders of magnitude in less than 20 V-cycles, although the multilevel convergence
is not exactly uniform with point density. The polynomial degree for restriction/prolongation and
discretization on the coarse levels is kept at three. Increasing this did not improve the
convergence shown in \cref{Fig:Dirichlet_Square}. Alternative cycles such as the W-cycle and adaptive V-cycle may
improve convergence but have not been explored in current study.

\begin{figure}[H]
	\centering
	\begin{subfigure}[t]{0.32\textwidth}
		\includegraphics[width=\textwidth]{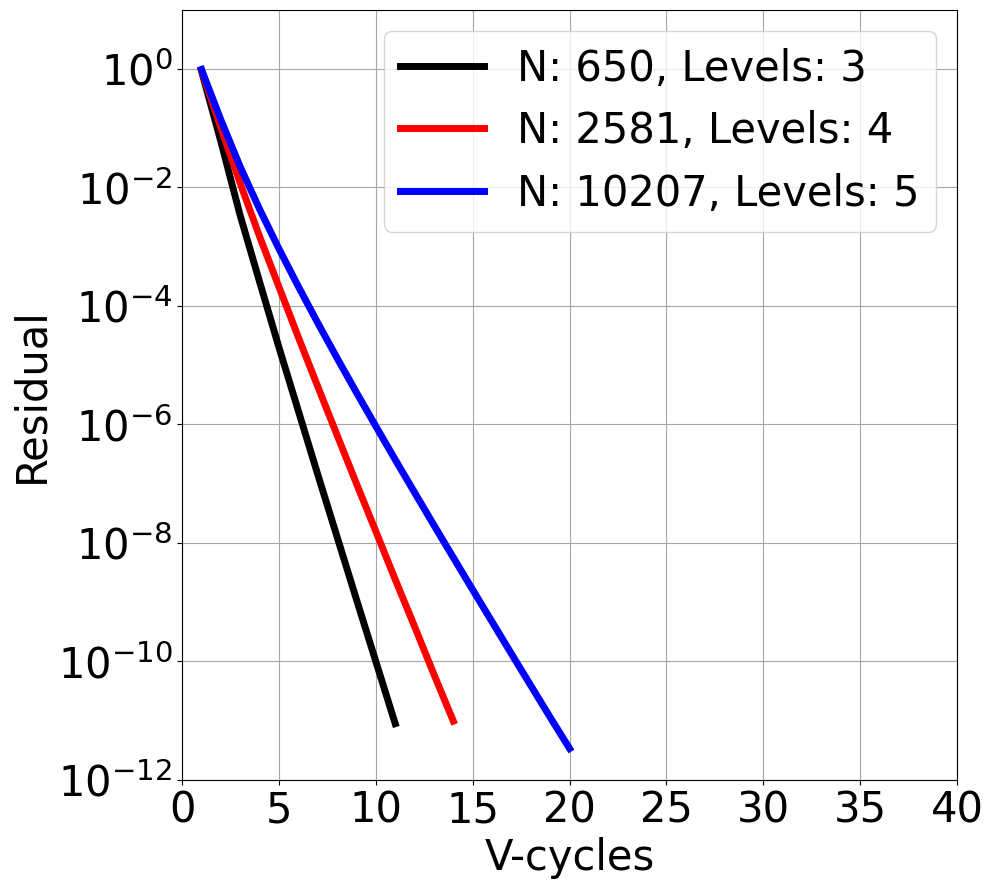}
		\caption{Degree of appended polynomial = 3}
	\end{subfigure}
	\begin{subfigure}[t]{0.32\textwidth}
		\includegraphics[width=\textwidth]{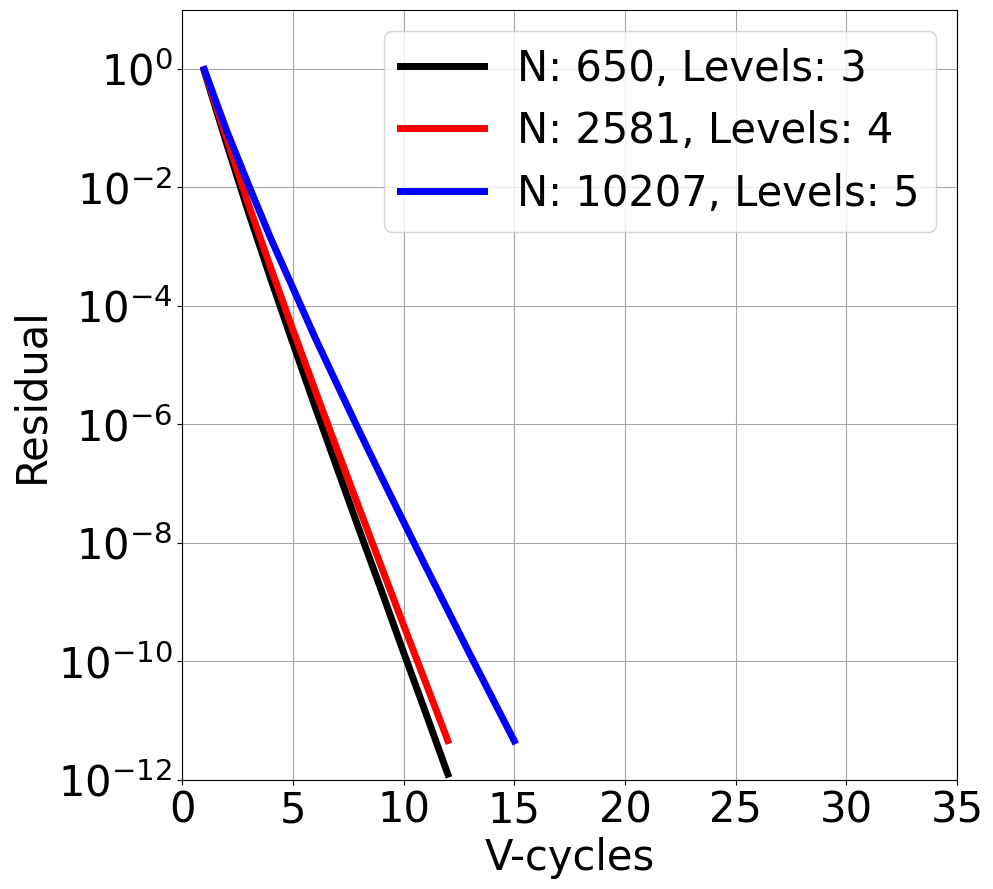}
		\caption{Degree of appended polynomial = 4}
	\end{subfigure}
	\begin{subfigure}[t]{0.32\textwidth}
		\includegraphics[width=\textwidth]{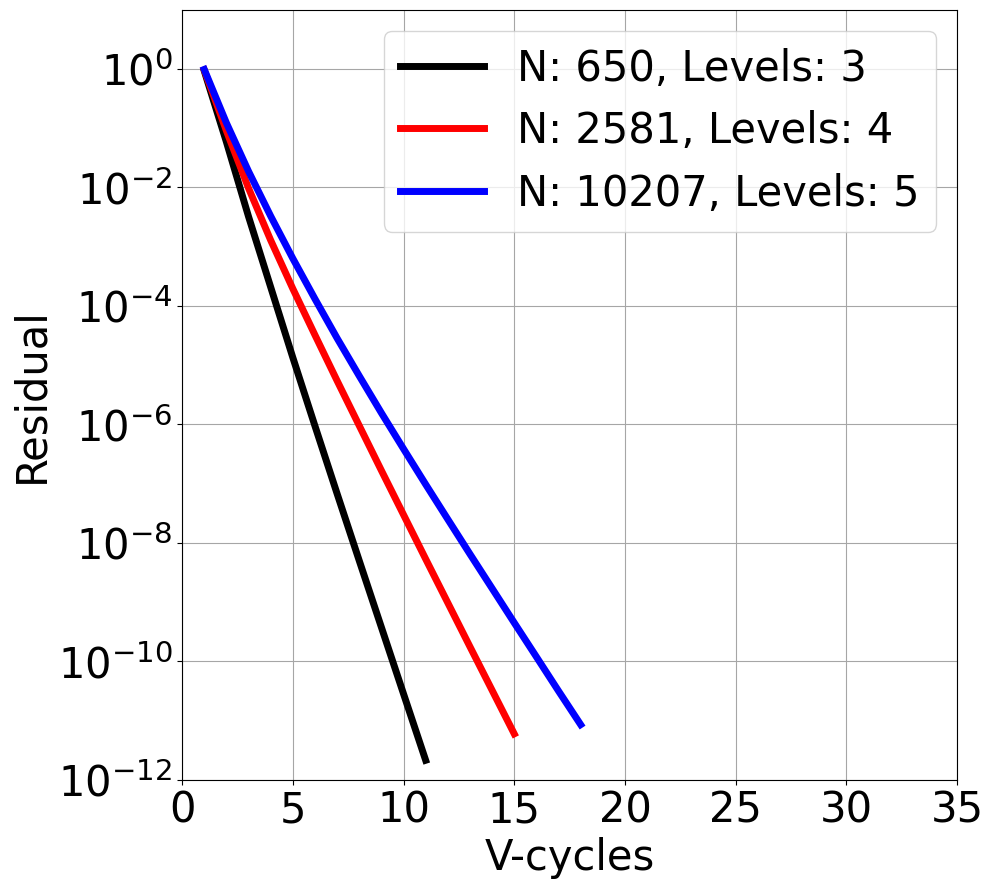}
		\caption{Degree of appended polynomial = 6}
	\end{subfigure}
	\caption{Convergence of the residual on a concentric annulus with Dirichlet boundary conditions}
	\label{Fig:Dirichlet_Conc_Circle}
\end{figure}

\begin{figure}[H]
	\centering
	\begin{subfigure}[t]{0.32\textwidth}
		\includegraphics[width=\textwidth]{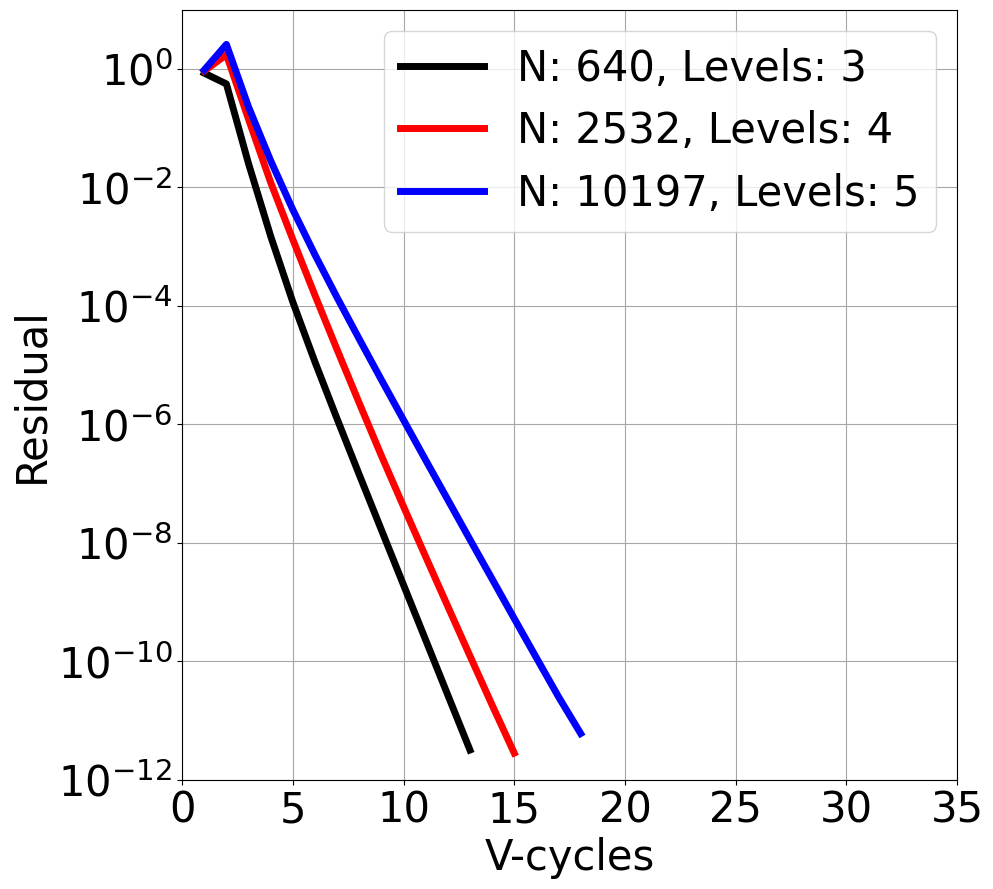}
		\caption{Degree of appended polynomial = 3}
	\end{subfigure}
	\begin{subfigure}[t]{0.32\textwidth}
		\includegraphics[width=\textwidth]{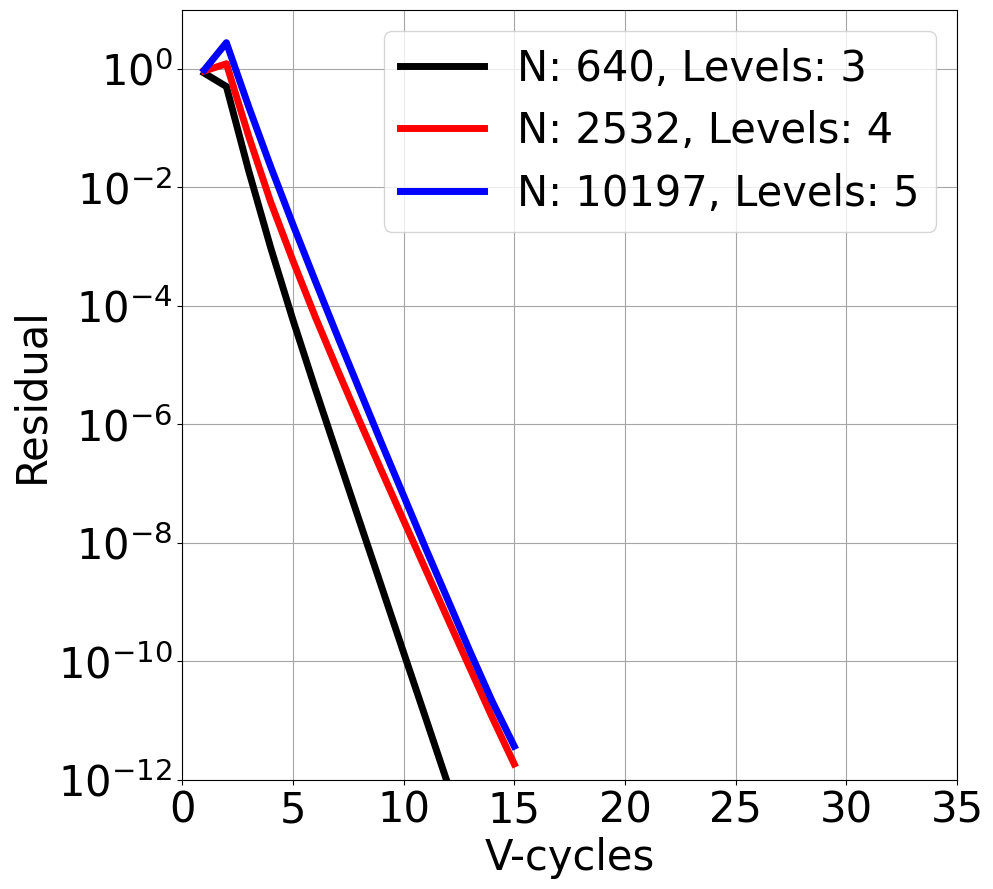}
		\caption{Degree of appended polynomial = 4}
	\end{subfigure}
	\begin{subfigure}[t]{0.32\textwidth}
		\includegraphics[width=\textwidth]{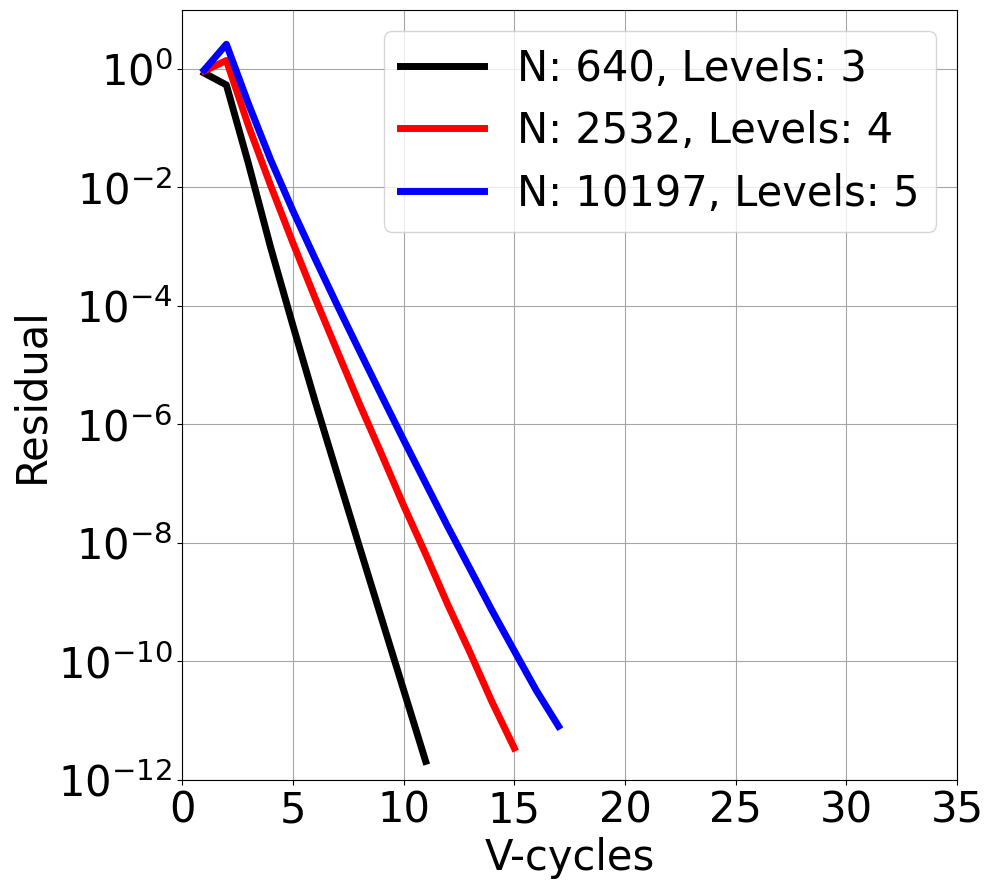}
		\caption{Degree of appended polynomial = 6}
	\end{subfigure}
	\caption{Convergence of the residual on a square-with-hole geometry with Dirichlet boundary conditions}
	\label{Fig:Dirichlet_Square_Hole}
\end{figure}

\Cref{Fig:Dirichlet_Conc_Circle,Fig:Dirichlet_Square_Hole} show similar convergence histories for the concentric annulus and the square-with-hole geometries respectively. The rates of convergence for these problems are also fast as
convergence is achieved in typically 10 to 20 V-cycles for ten orders of residual reduction.
However, convergence for larger point sets required more V cycles than for smaller sets,
displaying some multigrid inefficiency. The coefficients in the meshless discretization are
oscillatory, and the SOR with 5 relaxation sweeps is probably not sufficient to efficiently
annihilate the high frequency errors and may require more relaxation sweeps.

\begin{figure}[H]
	\centering
	\begin{subfigure}[t]{0.32\textwidth}
		\includegraphics[width=\textwidth]{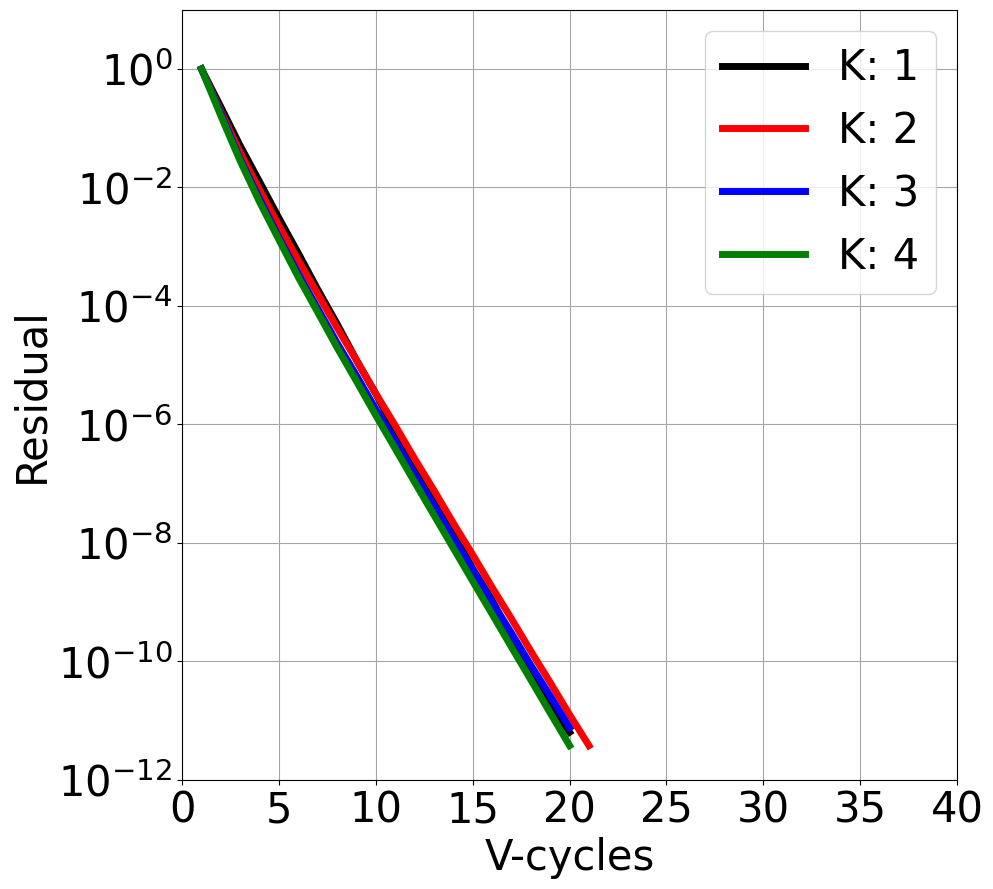}
		\caption{Square, Degree of appended polynomial = 6}
	\end{subfigure}
	\begin{subfigure}[t]{0.32\textwidth}
		\includegraphics[width=\textwidth]{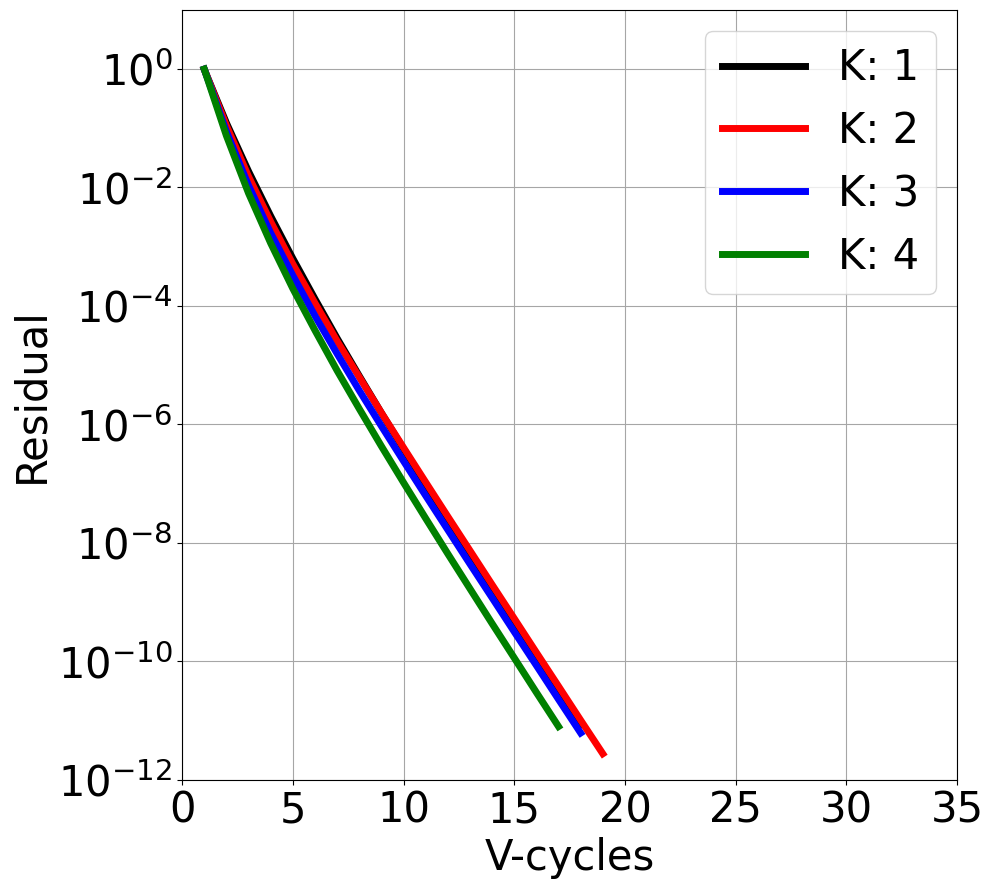}
		\caption{Concentric annulus, Degree of appended polynomial = 6}
	\end{subfigure}
	\begin{subfigure}[t]{0.32\textwidth}
		\includegraphics[width=\textwidth]{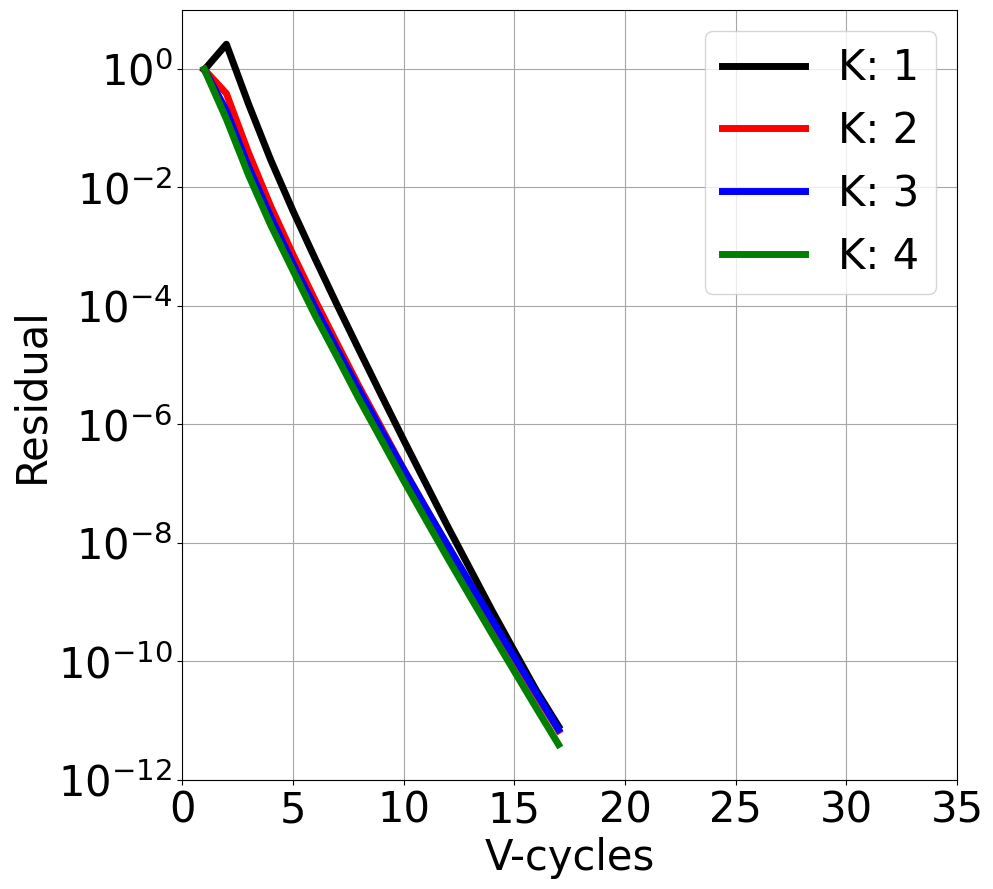}
		\caption{Square-with-hole, Degree of appended polynomial = 6}
	\end{subfigure}
	\caption{Effect of the wavenumber of the manufactured solution on multilevel convergence for finest level}
	\label{Fig:Frequency_Effect}
\end{figure}

\Cref{Fig:Frequency_Effect} shows the rates of convergence for different wavenumbers of the manufactured solution.
Higher wavenumbers have more oscillatory behavior and therefore can converge faster. The
results presented correspond to a point set consisting of around $\sim10,000$ points with an appended
polynomial of degree 6. The coarsest set has  $\sim90$ points (5 levels). We observe that all wavenumbers converge at nearly the same rate. Thus, henceforth, the calculations are performed by fixing the wavenumber $k$ at 1.

\begin{figure}[H]
	\centering
	\begin{subfigure}[t]{0.32\textwidth}
		\includegraphics[width=\textwidth]{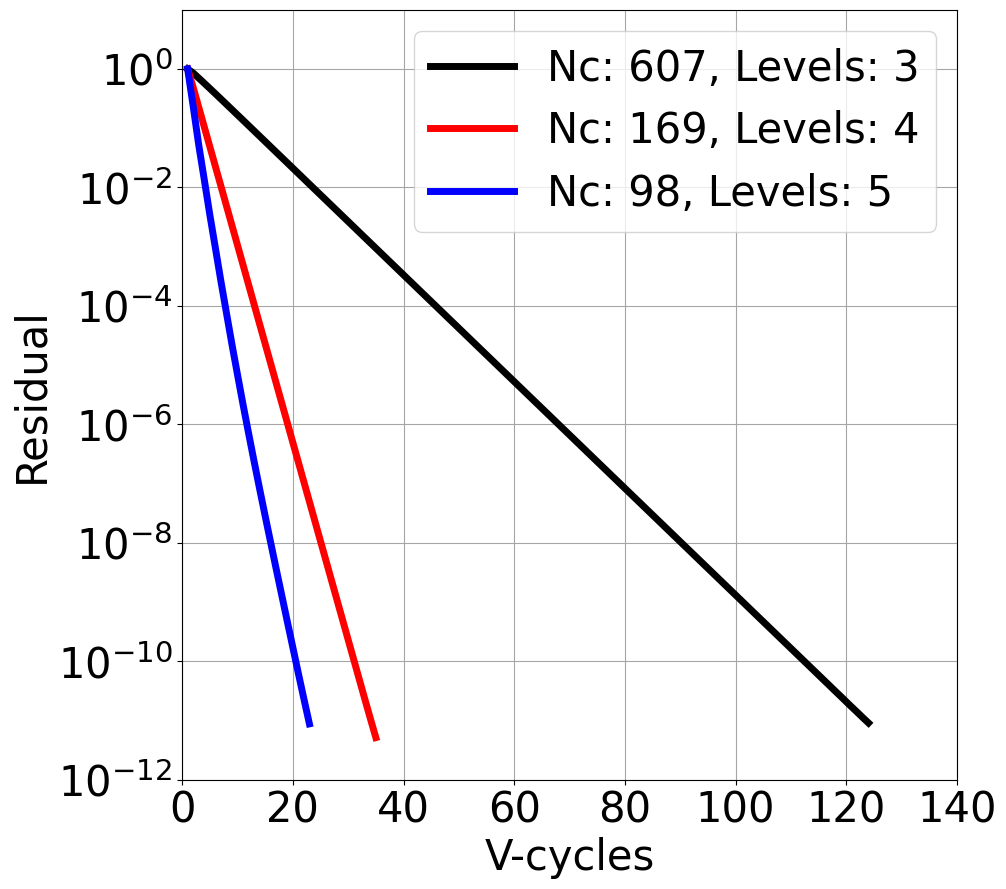}
	\end{subfigure}
	\caption{Effect of coarsest level on multilevel convergence}
	\label{Fig:Coarse_Effect}
\end{figure}

\Cref{Fig:Coarse_Effect} shows the rates of convergence for the square geometry with different degrees of
coarsening. The finest level (5) was kept the same with 10,023 points and the appended polynomial degree was fixed at $l$=3. The coarsest level
was varied from 1 (98 points) to 3 (607 points). It is seen that a factor of 5 reduction in the
number of V cycles is obtained by cycling up to level 1 versus cycling only to the third level. This
factor was slightly different for other problems, and for other degrees of appended polynomials.
As expected, convergence was better with the coarsest level having 98 points. The number of
iterations with only a single level relaxation was significantly greater and not included here.

\subsection{Results for All Neumann Boundary Conditions}

\begin{figure}[H]
	\centering
	\begin{subfigure}[t]{0.32\textwidth}
		\includegraphics[width=\textwidth]{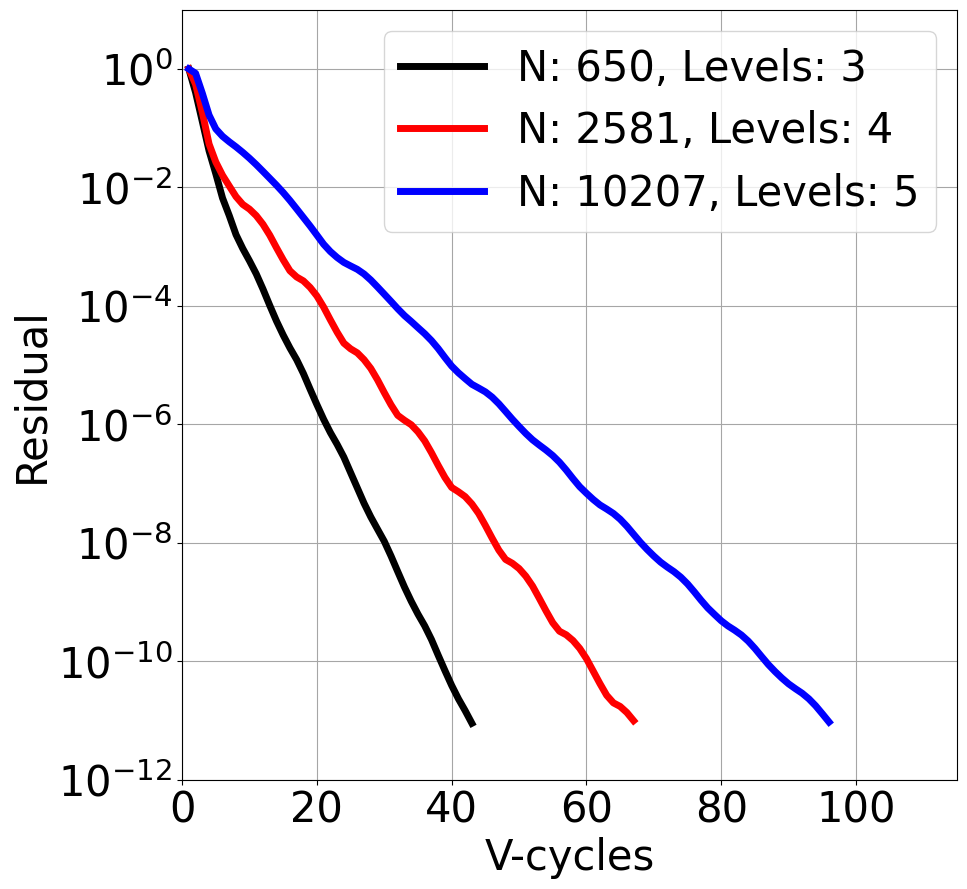}
		\caption{Degree of appended polynomial = 3}
	\end{subfigure}
	\begin{subfigure}[t]{0.32\textwidth}
		\includegraphics[width=\textwidth]{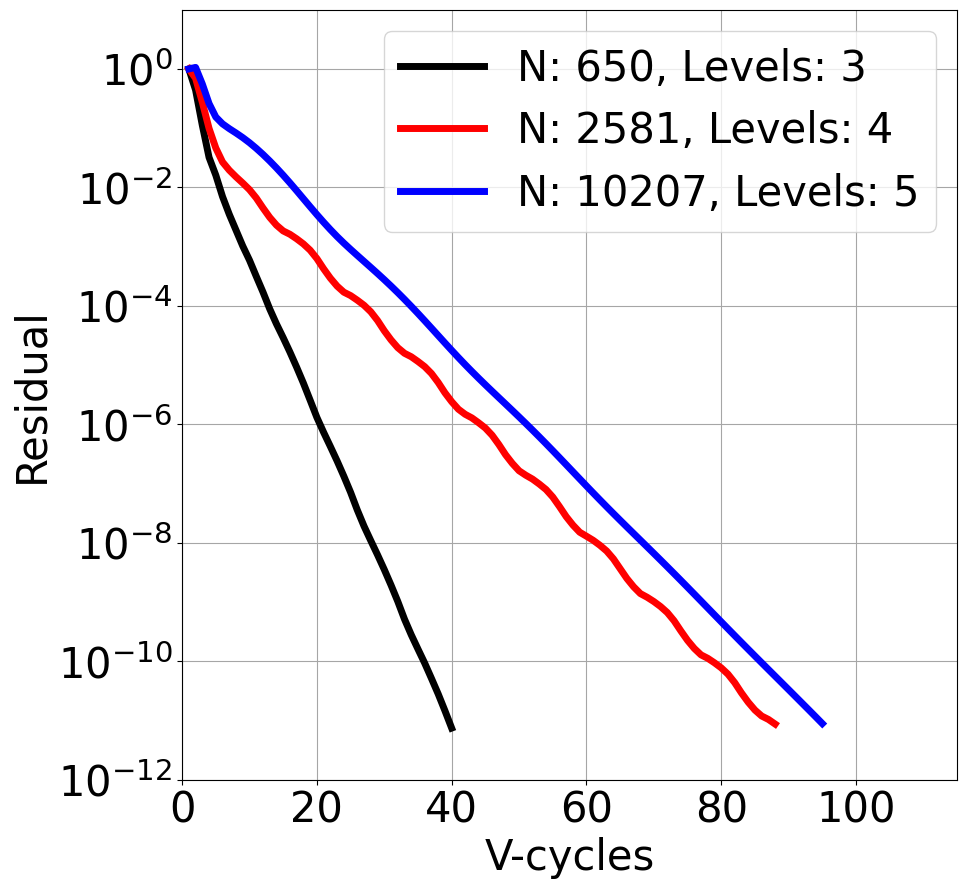}
		\caption{Degree of appended polynomial = 4}
	\end{subfigure}
	\begin{subfigure}[t]{0.32\textwidth}
		\includegraphics[width=\textwidth]{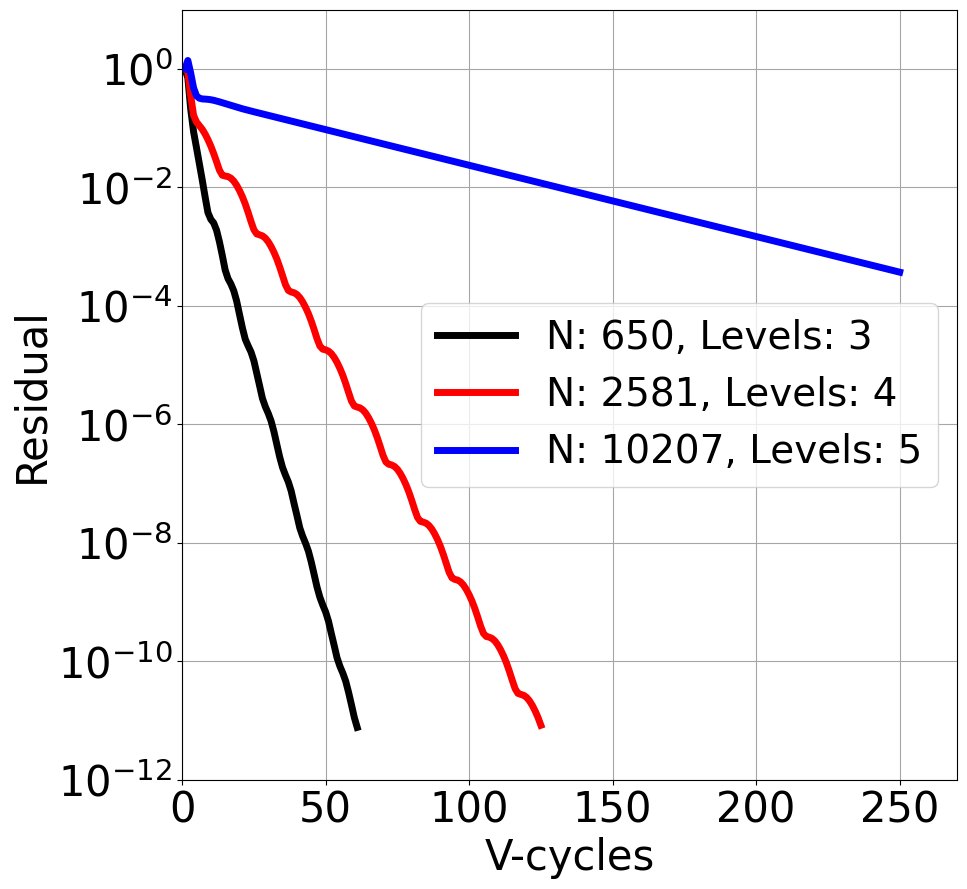}
		\caption{Degree of appended polynomial = 5}
	\end{subfigure}
	\caption{Convergence of the residual on a concentric annulus with all-Neumann boundary conditions}
	\label{Fig:Neumann_Conc_Circle}
\end{figure}

\begin{figure}[H]
	\centering
	\begin{subfigure}[t]{0.32\textwidth}
		\includegraphics[width=\textwidth]{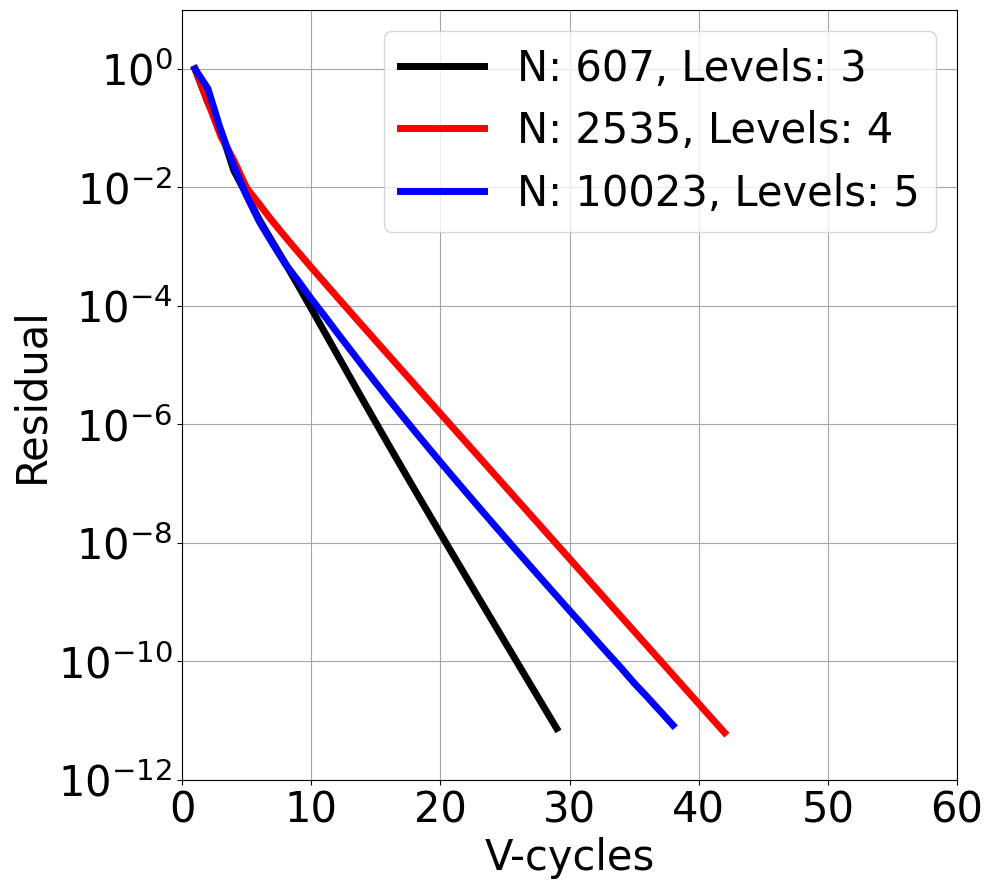}
		\caption{Degree of appended polynomial = 3}
	\end{subfigure}
	\begin{subfigure}[t]{0.32\textwidth}
		\includegraphics[width=\textwidth]{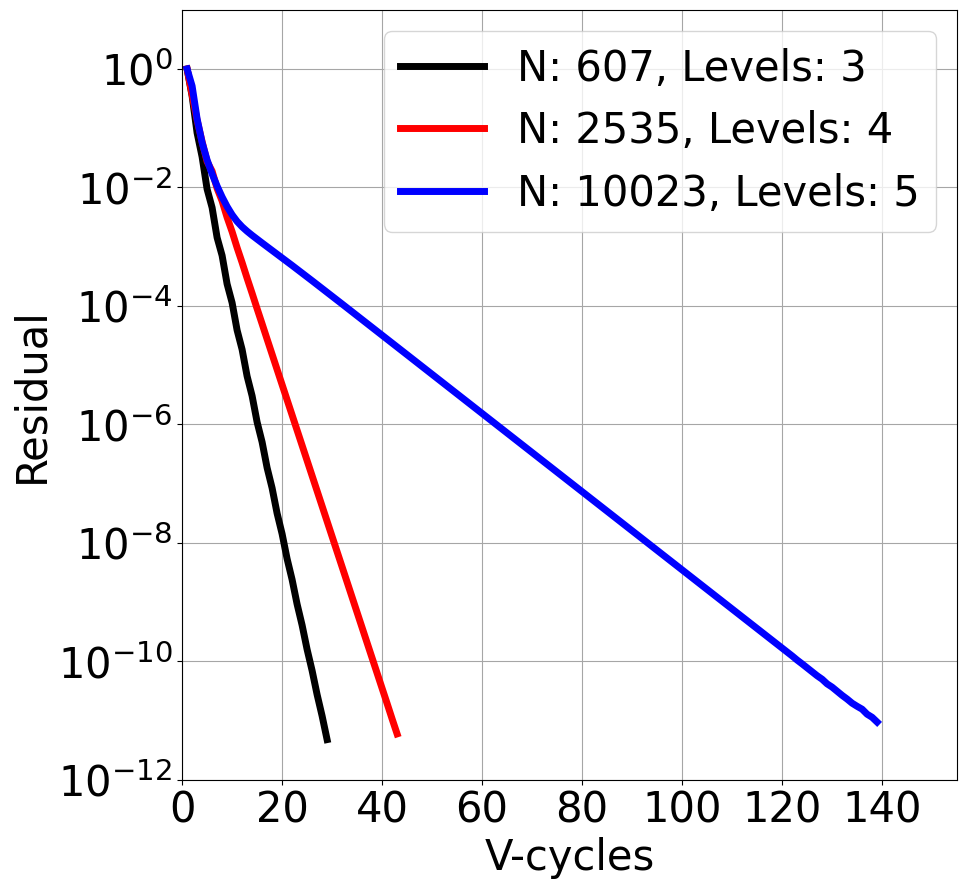}
		\caption{Degree of appended polynomial = 4}
	\end{subfigure}
	\begin{subfigure}[t]{0.32\textwidth}
		\includegraphics[width=\textwidth]{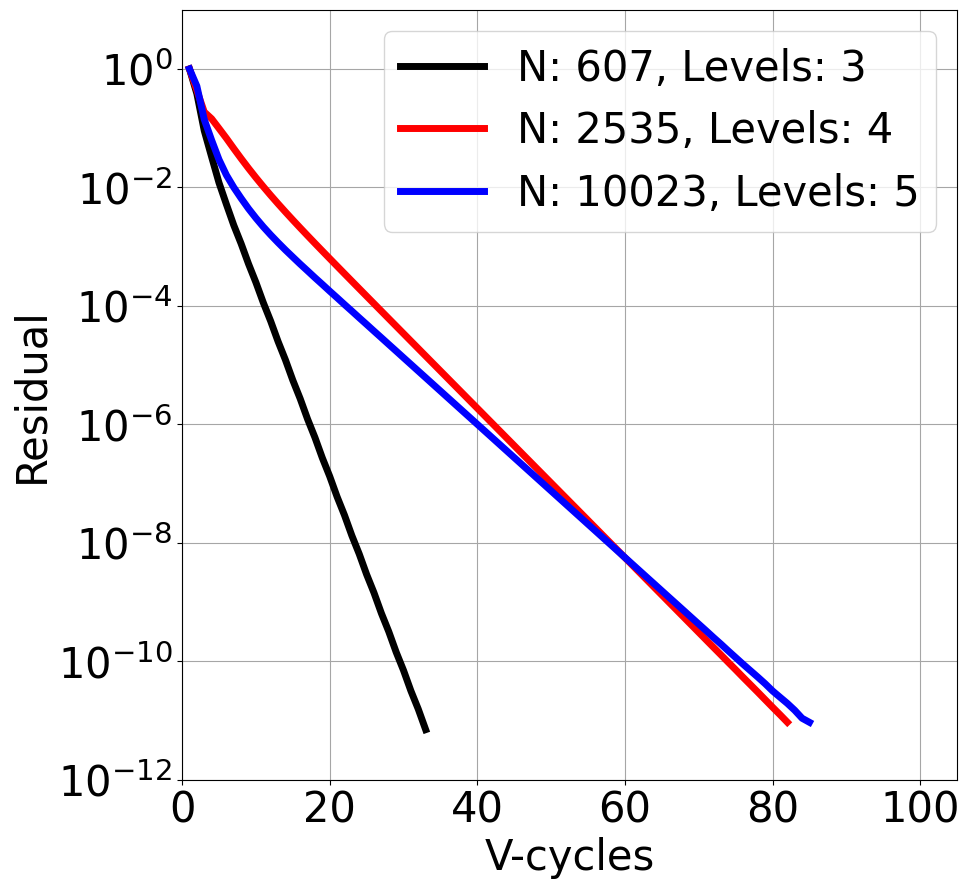}
		\caption{Degree of appended polynomial = 5}
	\end{subfigure}
	\caption{Convergence of the residual on a square geometry with all-Neumann boundary conditions}
	\label{Fig:Neumann_Square}
\end{figure}

\begin{figure}[H]
	\centering
	\begin{subfigure}[t]{0.32\textwidth}
		\includegraphics[width=\textwidth]{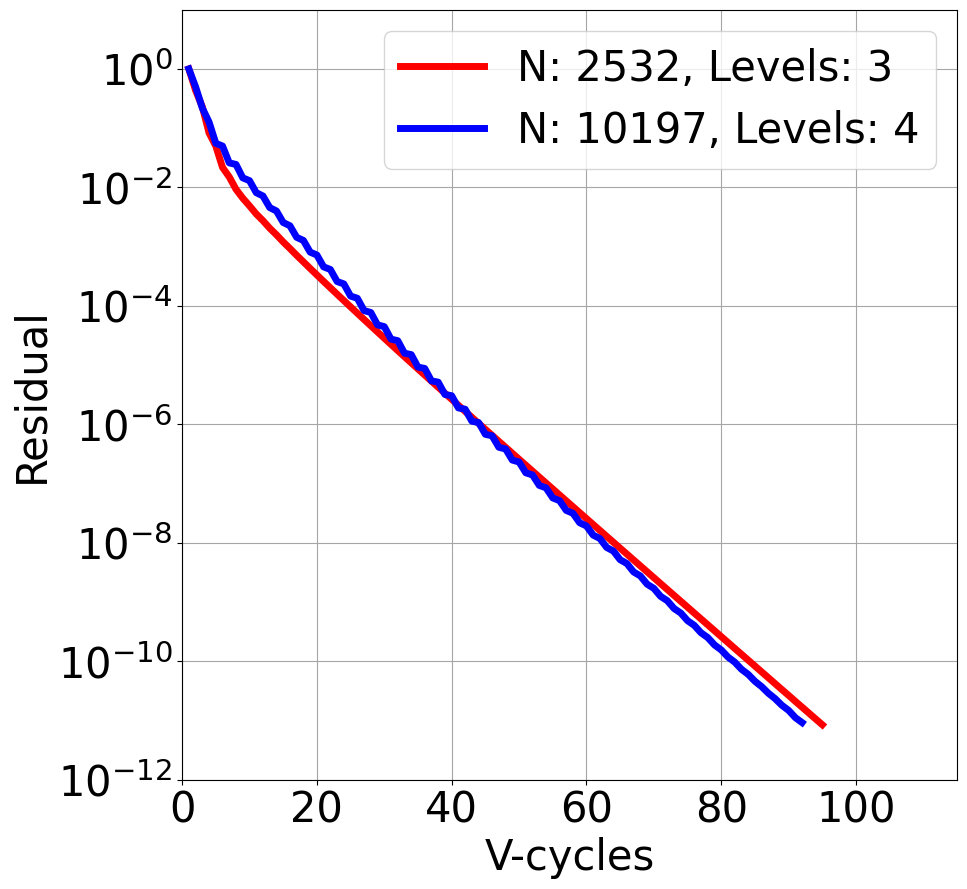}
		\caption{Degree of appended polynomial = 3}
	\end{subfigure}
	\begin{subfigure}[t]{0.32\textwidth}
		\includegraphics[width=\textwidth]{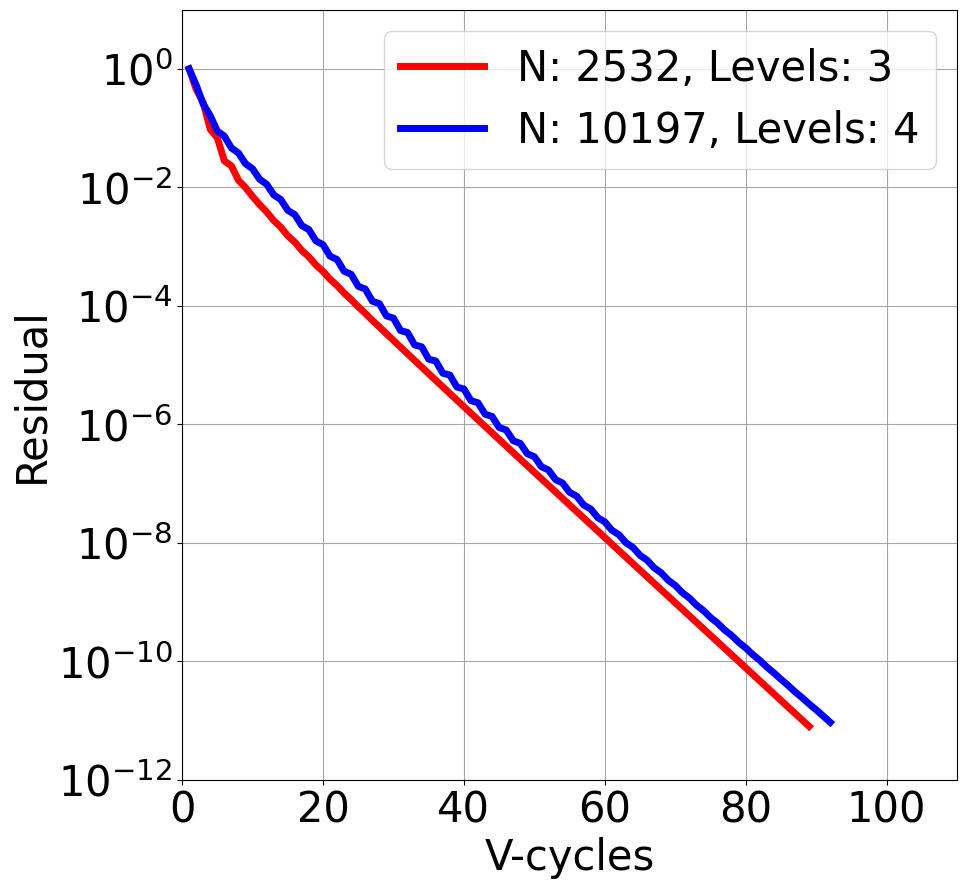}
		\caption{Degree of appended polynomial = 4}
	\end{subfigure}
	\begin{subfigure}[t]{0.32\textwidth}
		\includegraphics[width=\textwidth]{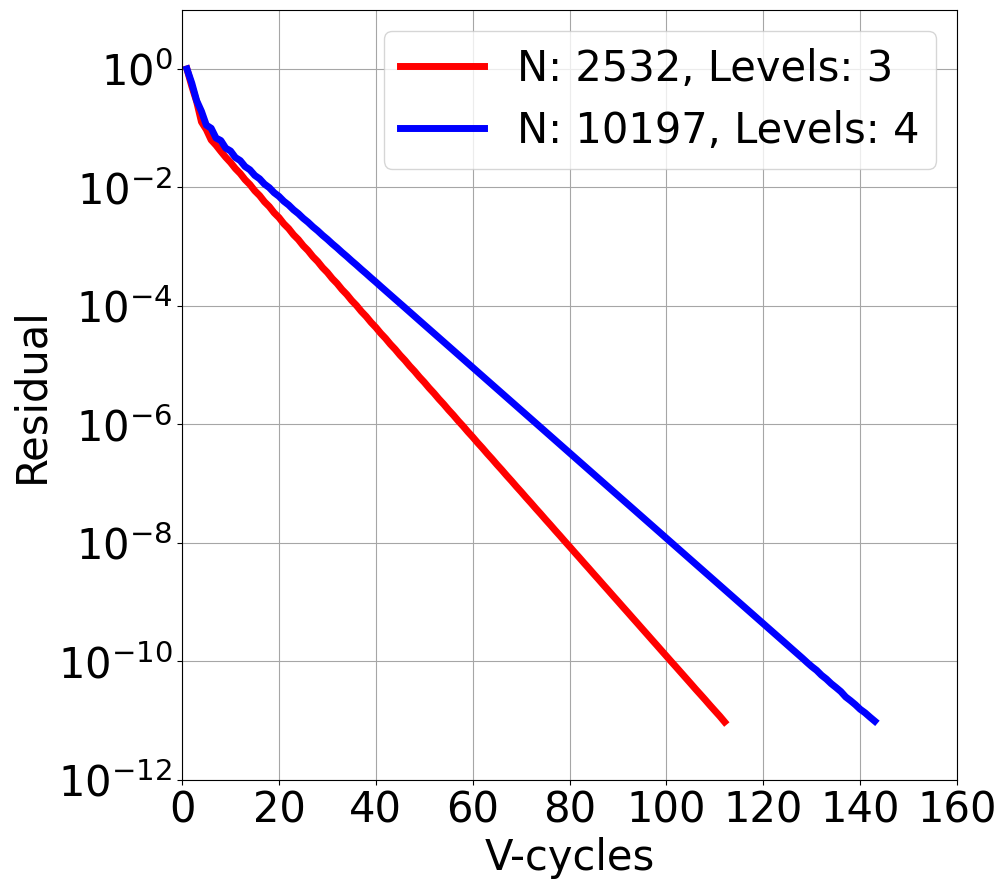}
		\caption{Degree of appended polynomial = 5}
	\end{subfigure}
	\caption{Convergence of the residual on a square-with-hole geometry with all-Neumann boundary conditions (coarsest level = 176 points)}
	\label{Fig:Neumann_Square_Hole}
\end{figure}

We next present results for Neumann conditions prescribed on all the boundaries. The analytical
derivative given by the manufactured solution was imposed at all the points on the boundaries.
Note that when all the boundaries are of Neumann type, the problem is ill-conditioned as the
solution can only be computed to an arbitrary constant. Hence the matrix must be regularized
either by fixing one of the points to some arbitrary value or imposing a global constraint on the
values. Here we use the regularization procedure in which the sum of all values was assigned an
arbitrary level of zero \cite{regularize2019}. The solution matrix was therefore modified by adding an extra
equation of unity coefficients and zero right-hand side. \Cref{Fig:Neumann_Conc_Circle,Fig:Neumann_Square,Fig:Neumann_Square_Hole}  show the multilevel
convergence for the three problems with all Neumann boundary conditions. SOR
with over-relaxation factor of 1.4 was used with five (5) iterations. It is seen that even after
regularization, the convergence rate is not satisfactory. Compared with the case of Dirichlet
boundary condition, the Neumann condition took more V-cycles to converge, and the
convergence depended on the refinement. Some of the computations even diverged, especially
for high degree of polynomial ($l$ = 6). In the case of square-with-hole geometry, limiting the coarsest level to 176 points gave better convergence.

\subsection{Multilevel Preconditioned GMRES}

The all-Neumann boundary condition is important in many fields including for determination of
the pressure field in simulations of incompressible fluid flows with prescribed boundary conditions. Since the multilevel solver with
pure SOR relaxation converged slowly for cases of the all-Neumann boundary conditions, we explored a
combination of GMRES and the multilevel algorithm in which the multilevel algorithm is used to
resolve the residuals computed inside a Krylov Subspace Projector (KSP). Five sweeps of the SOR are also performed at each level, but residuals to be smoothed are calculated by the GMRES
algorithm on the finest level. The pseudo code for this algorithm is given below:\\

\noindent \textbf{Pseudocode for Multilevel Preconditioned GMRES Algorithm }

\noindent Let $x = 0$ , $r = b$, $k = 0$, tol = $10^{-10}$ \\
 \noindent while(true)
    \par $k++$
    \par Solve $M_{g} z = r$
    \par if($k == 1$)
         \par \hspace{1cm} $ p_{k} = z$
          \par \hspace{1cm} $w_{k} = A p_{k}$
    \par else
       \par \hspace{1cm} $p_{k} = z$
         \par \hspace{1cm} for $i = 1...k-1$
          \par \hspace{2cm} $ p_{k} = p_{k} – [(w_{i})^{T} (A z) / (w_{i})^{T}(w_{i})] $
         \par \hspace{1cm} $w_{k} = A p_{k}$
    \par$\alpha = (w_{k})^{T} (Ar) / (w_{k})^{T} (w_{k})$
    \par$x = x + \alpha p_{k}$
    \par$r = r – \alpha w_{k}$
    \par if ( $(||r|| / ||b||) < $  tol )
       \par \hspace{1cm} break \\

\par Here, $p_{k}$’s are mutually orthogonal in the $A^{T} A$ norm, which is Symmetric Positive Definite (SPD)
even if $A$ is not SPD. Thus, the updated solution is the best fit in the space of the $p_{k}$’s in the $A^{T}A$
norm. Note that we have used modified Gram-Schmidt procedure to orthogonalize the search
directions for better numerical stability. The solution to $M_{g} z = r$ applies one V-cycle of the
multilevel operator $M_{g}$ with the righthand side vector as $r$. The multilevel procedure starts with
an initial vector of zero corrections and stores the result after one V-cycle in the vector $z$. This
procedure was seen to be efficient and robust for all the pure Neumann boundary condition
cases. Also, the number of V-cycles was not observed to increase very much with polynomial degree, although the condition number increases significantly with the degree of the appended polynomial.

\begin{figure}[H]
	\centering
	\begin{subfigure}[t]{0.32\textwidth}
		\includegraphics[width=\textwidth]{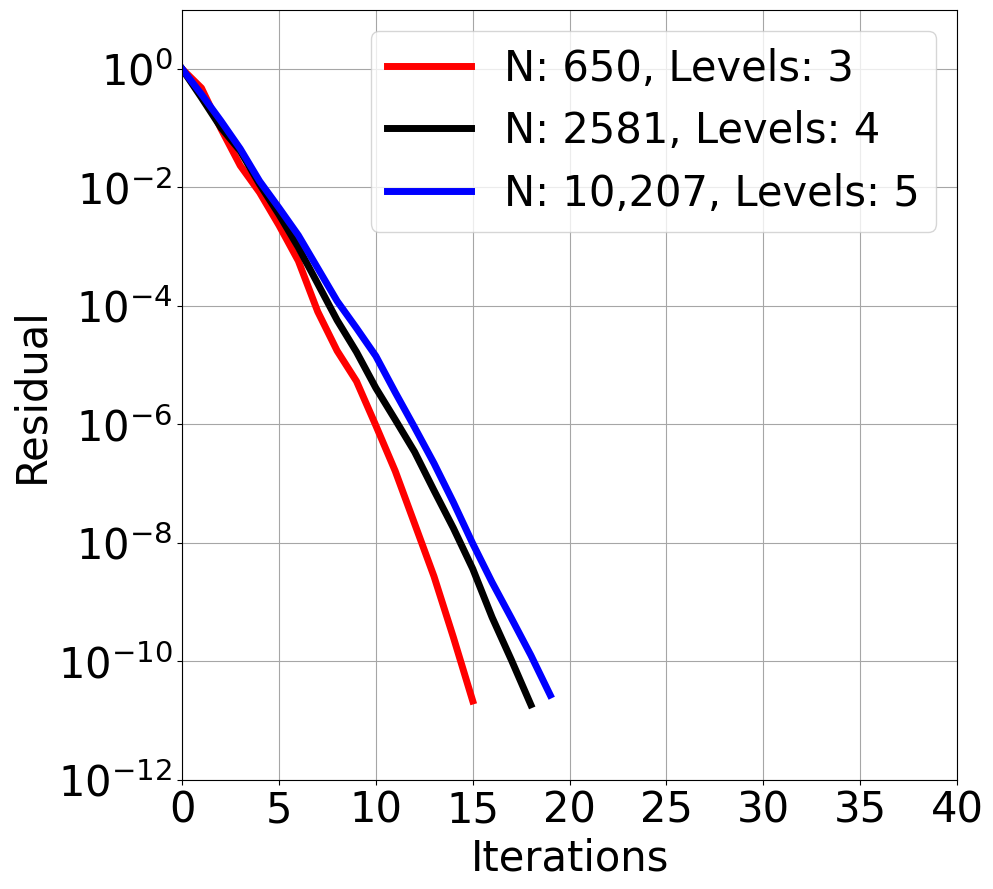}
		\caption{Degree of appended polynomial = 4}
	\end{subfigure}
	\begin{subfigure}[t]{0.32\textwidth}
		\includegraphics[width=\textwidth]{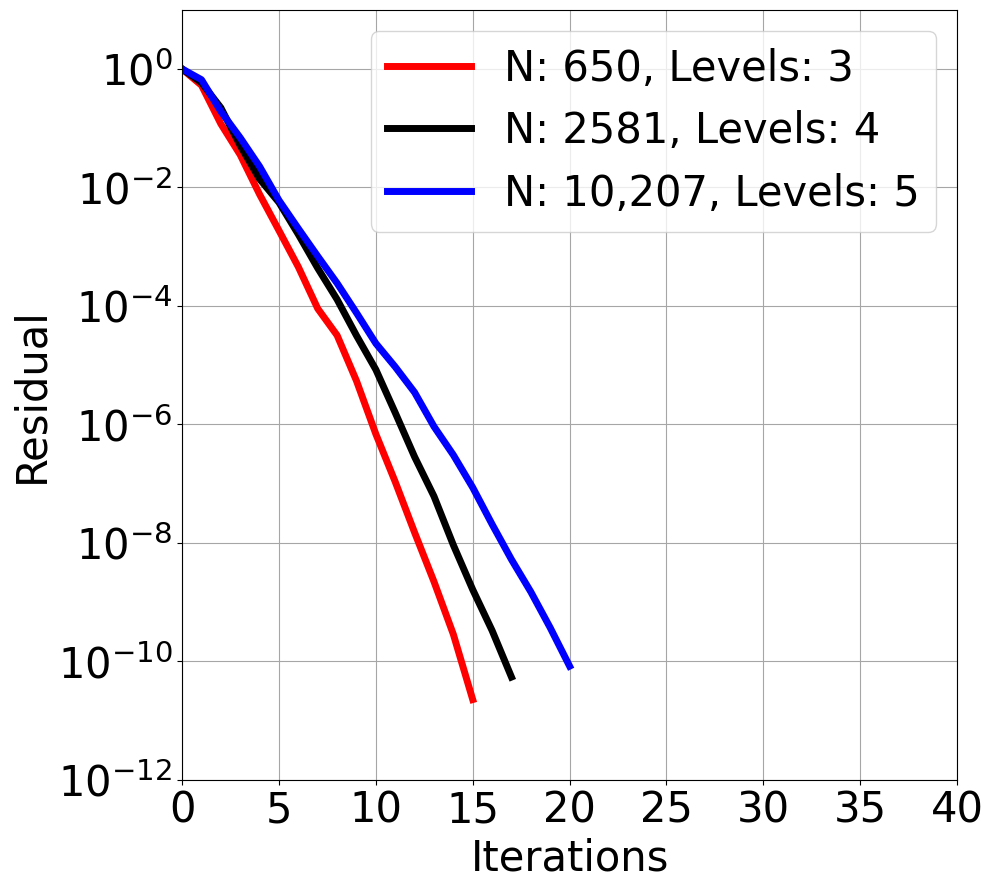}
		\caption{Degree of appended polynomial = 5}
	\end{subfigure}
	\begin{subfigure}[t]{0.32\textwidth}
		\includegraphics[width=\textwidth]{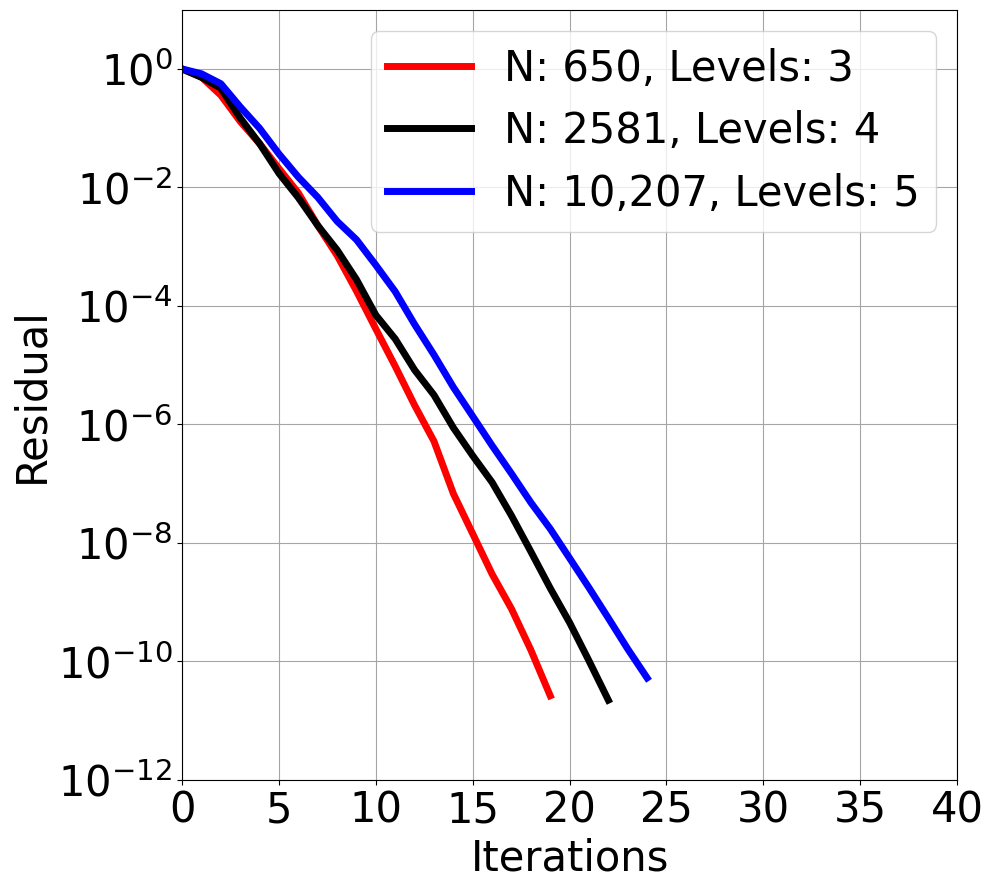}
		\caption{Degree of appended polynomial = 6}
	\end{subfigure}
	\caption{Convergence of the GMRES-Multilevel algorithm for the concentric annulus with all-Neumann boundary conditions}
	\label{Fig:Neumann_Conc_Circle_GMRES}
\end{figure}

\begin{figure}[H]
	\centering
	\begin{subfigure}[t]{0.32\textwidth}
		\includegraphics[width=\textwidth]{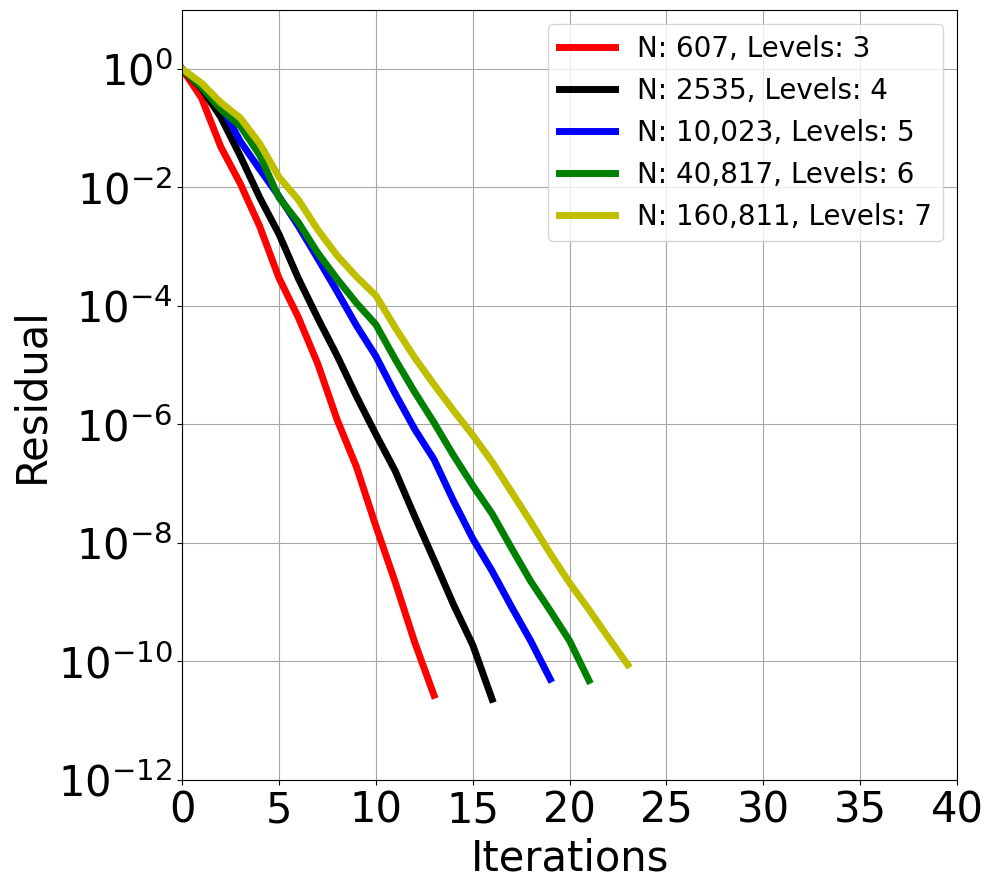}
		\caption{Degree of appended polynomial = 4}
	\end{subfigure}
	\begin{subfigure}[t]{0.32\textwidth}
		\includegraphics[width=\textwidth]{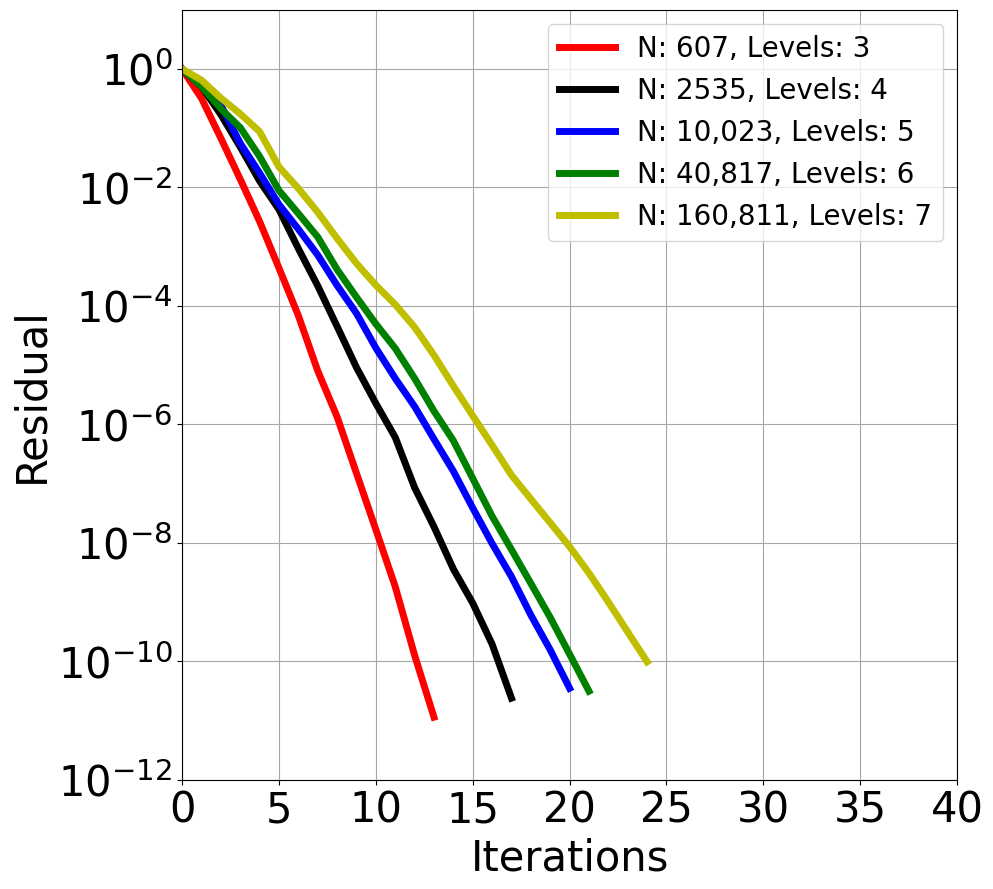}
		\caption{Degree of appended polynomial = 5}
	\end{subfigure}
	\begin{subfigure}[t]{0.32\textwidth}
		\includegraphics[width=\textwidth]{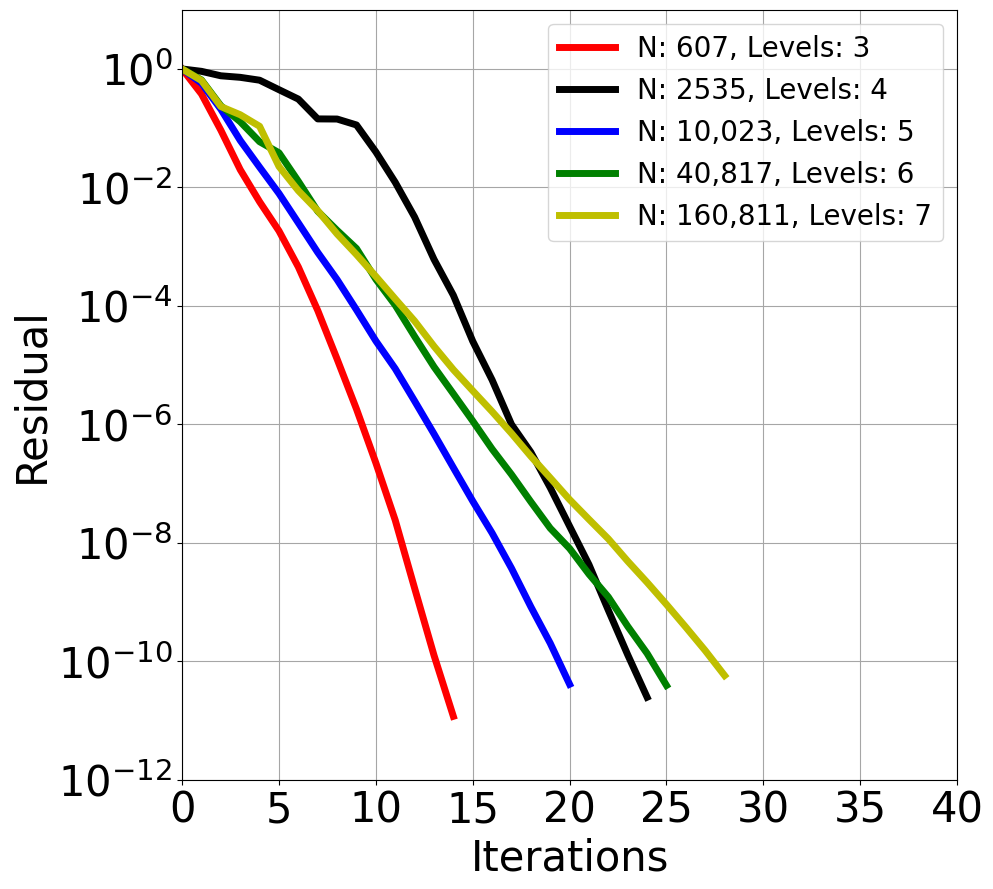}
		\caption{Degree of appended polynomial = 6}
	\end{subfigure}
	\caption{Convergence of the GMRES-Multilevel algorithm for the square geometry with all-Neumann boundary conditions}
	\label{Fig:Neumann_Square_GMRES}
\end{figure}

\begin{figure}[H]
	\centering
	\begin{subfigure}[t]{0.32\textwidth}
		\includegraphics[width=\textwidth]{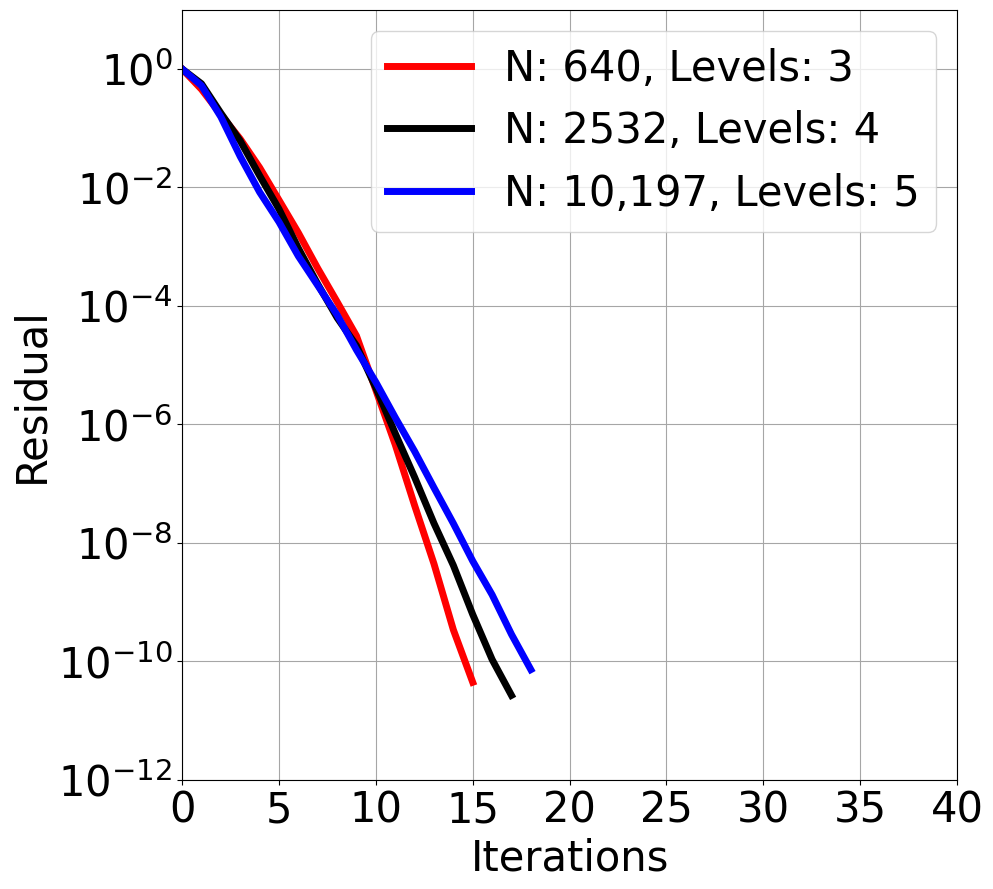}
		\caption{Degree of appended polynomial = 4}
	\end{subfigure}
	\begin{subfigure}[t]{0.32\textwidth}
		\includegraphics[width=\textwidth]{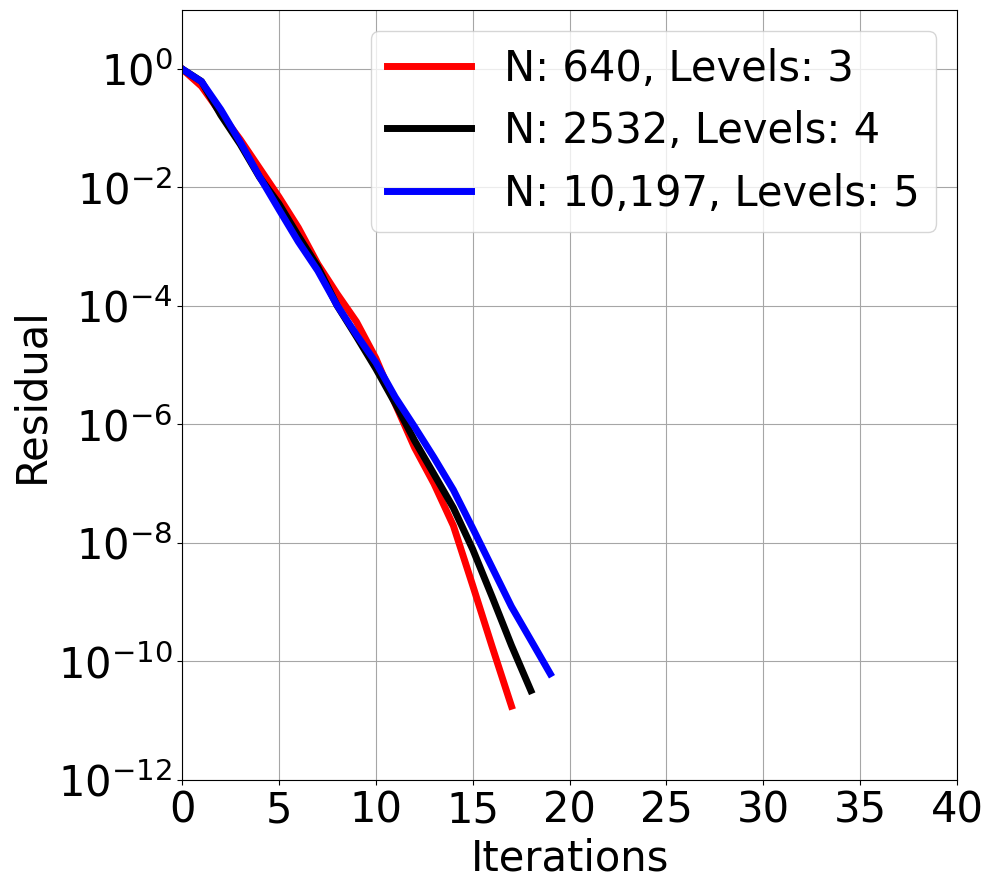}
		\caption{Degree of appended polynomial = 5}
	\end{subfigure}
	\begin{subfigure}[t]{0.32\textwidth}
		\includegraphics[width=\textwidth]{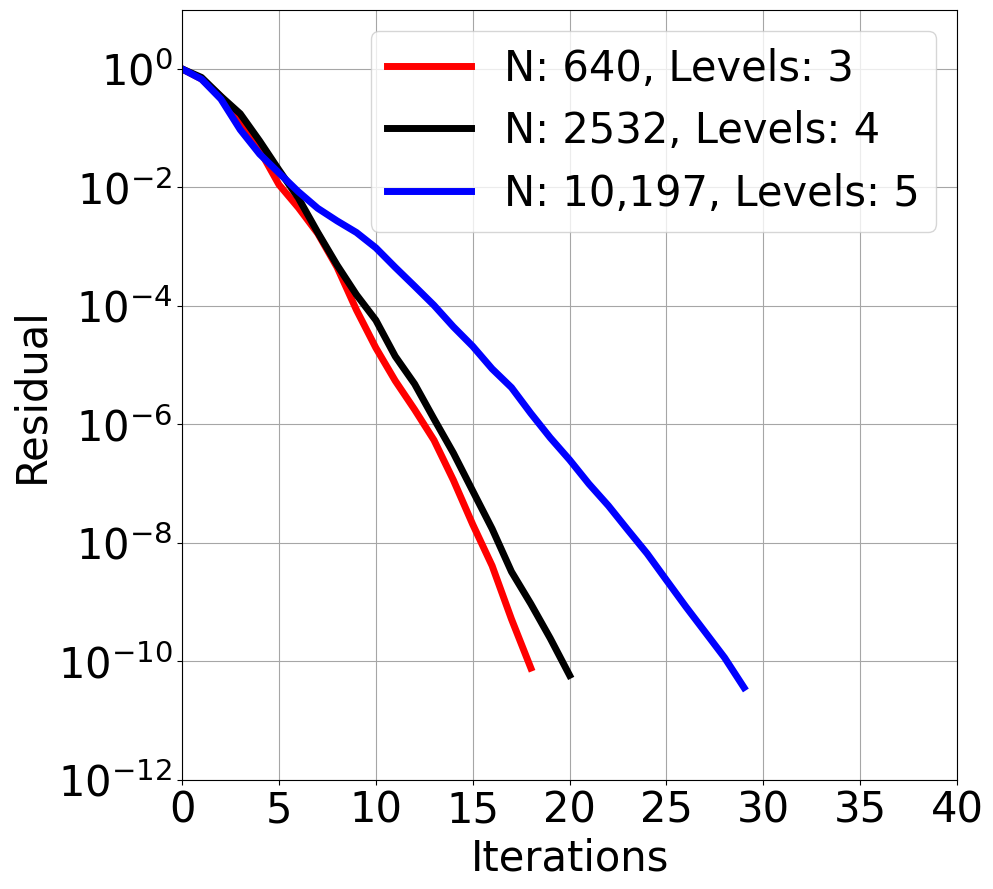}
		\caption{Degree of appended polynomial = 6}
	\end{subfigure}
	\caption{Convergence of the GMRES-Multilevel algorithm for the square-with-hole geometry with all-Neumann boundary conditions}
	\label{Fig:Neumann_Square_Hole_GMRES}
\end{figure}

\Cref{Fig:Neumann_Conc_Circle_GMRES,Fig:Neumann_Square_GMRES,Fig:Neumann_Square_Hole_GMRES} show convergence histories of this algorithm for the wavenumber $k$ = 1 and for
the three geometries. It is seen that the combination of GMRES with the multilevel solver for the
residuals gives fast convergence for the all-Neumann boundary condition. Although all point
sets do not converge at same rates, the multilevel property is nearly achieved. The slightly slower
convergence for higher order polynomial cases may be a result of the condition number, and the
large bandwidth of the sparse coefficient matrix. Other relaxation schemes such a block-SOR
may provide faster convergence and will be explored in future studies.

\subsection{Discretization Error}
As mentioned earlier, the PHS-RBF method displays exponential convergence with point spacing
as per the degree of the appended polynomial. The accuracy of the present solutions for the three
geometries is demonstrated by plotting below the differences from the known manufactured
solutions. For a given degree of the polynomial $l$, the second derivative displays at least ($l$-1)
degree of accuracy.

\begin{figure}[H]
	\centering
	\begin{subfigure}[t]{0.32\textwidth}
		\includegraphics[width=\textwidth]{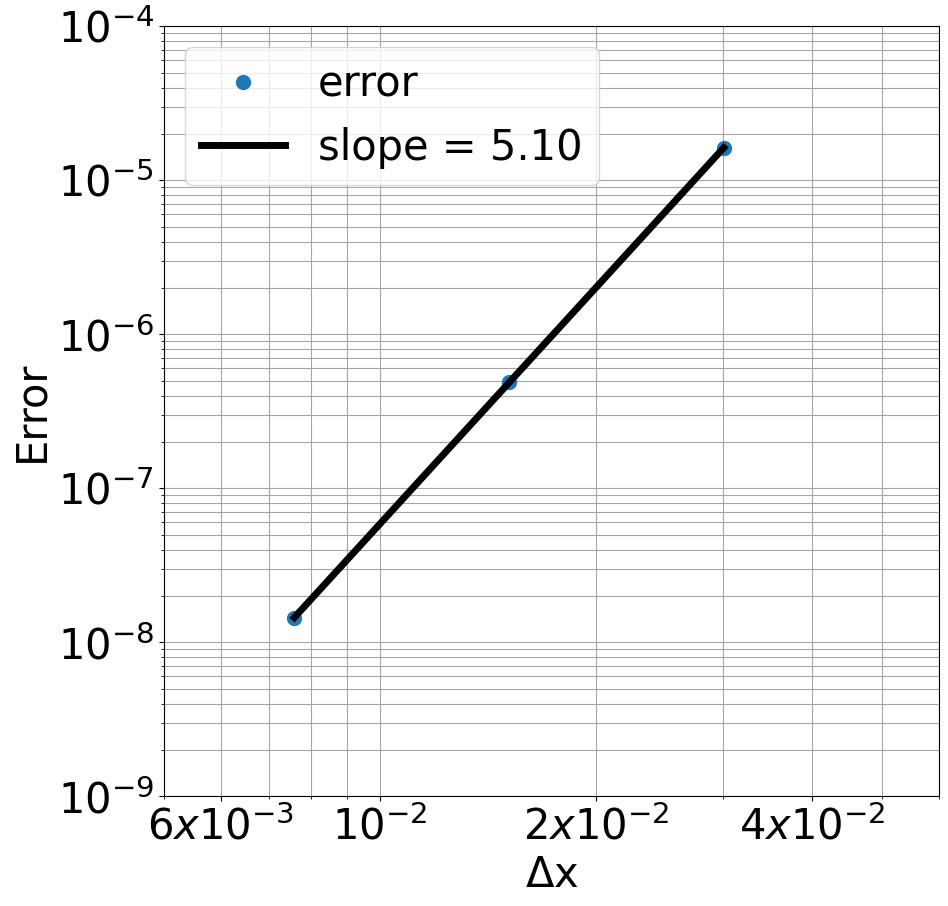}
		\caption{Concentric annulus}
	\end{subfigure}
	\begin{subfigure}[t]{0.32\textwidth}
		\includegraphics[width=\textwidth]{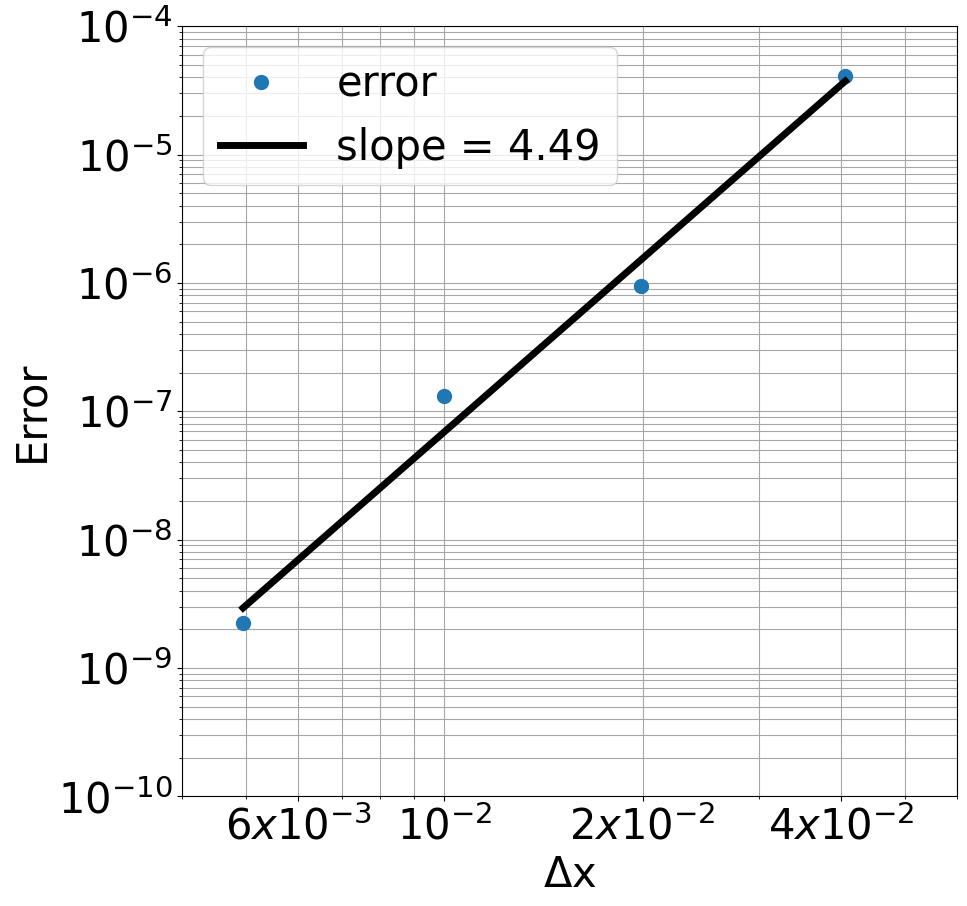}
		\caption{Square}
	\end{subfigure}
	\begin{subfigure}[t]{0.32\textwidth}
		\includegraphics[width=\textwidth]{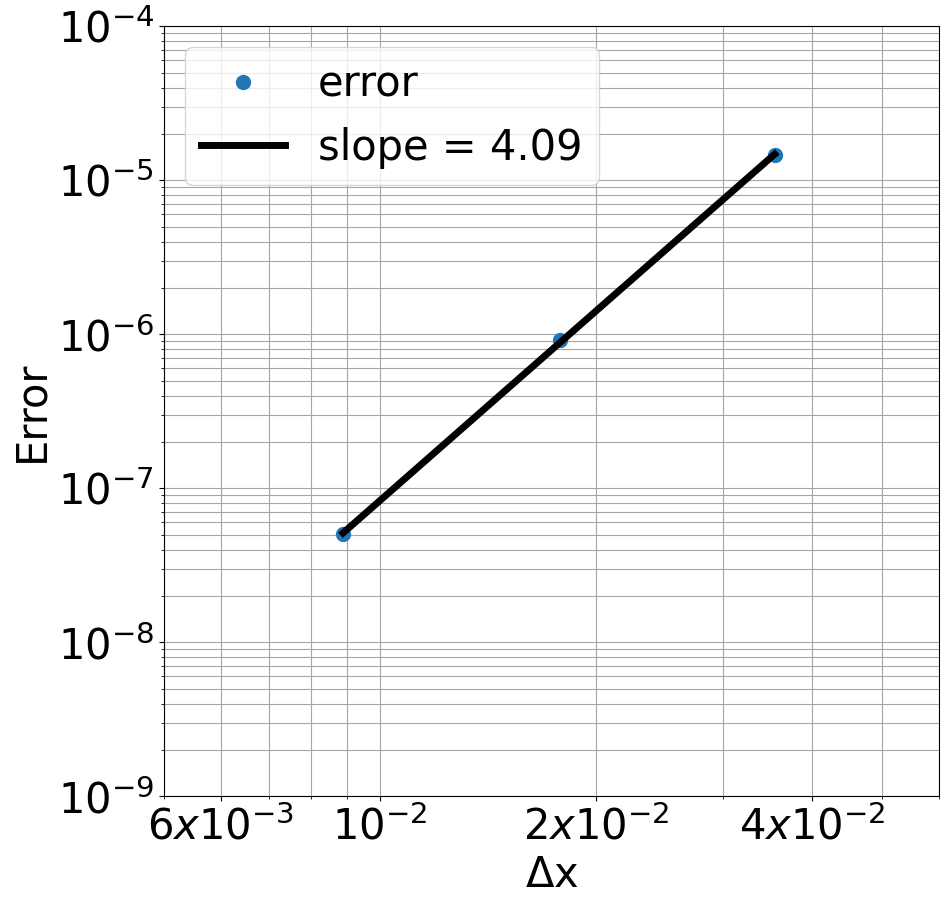}
		\caption{Square-with-hole}
	\end{subfigure}
	\caption{Discretization errors for the three geometries with degree of appended polynomial of 5.}
	\label{Fig:Error_Plots}
\end{figure}

\noindent The calculated errors for each case are plotted in \cref{Fig:Error_Plots} for polynomial degree of 5 against an
average spacing $\Delta x$ defined as:

\begin{equation}
    \Delta x = (Ar/n)^{0.5}
    \label{Eq:Spacing}
\end{equation}

\noindent where $Ar$ and $n$ in \cref{Eq:Spacing} refer to the area of the domain and number of points respectively.

\par The errors are
seen to follow the expected trends. It must be pointed out that both the point spacings and
polynomial order can be also varied locally as well as globally. This gives more freedom in targeting accuracy to desired regions through combined refinements in point spacing as well as
the polynomial degree.

\section{Summary and Future Research}
 We have presented in this paper a multilevel algorithm for solving the discrete equations
resulting from meshless discretization of the Poisson equation in a complex domain. The
meshless discretization is done by scattered data interpolation using polyharmonic splines with
appended polynomials of high degree. The multilevel algorithm uses non-nested point sets and restriction/prolongations using the PHS-RBF interpolations. Systematic tests have been
performed for three different geometries using manufactured solutions. The algorithm is seen
to converge rapidly for Dirichlet boundary conditions for large number of scattered points and
for high order appended polynomials. However, for the case of all-Neumann boundary condition,
the convergence with pure SOR is seen to be inferior and not be robust. The multilevel method
was then investigated as a preconditioner to a GMRES algorithm. Good and robust convergence
has been achieved for the all-Neumann cases. The Multilevel preconditioned GMRES algorithm is then implemented in a fractional step algorithm
to solve the pressure Poisson equation. Good convergence is seen for two model problems tested with approximately 10,000 scattered points (Kovasznay flow \cite{kovasznay1948laminar} and flow between two concentric cylinders due to a rotating inner cylinder). Future efforts will be directed towards extensions to three-dimensions
and parallelization of the entire algorithm including the coefficient computations.
Implementation on massively parallel architectures such as the GPU will also be pursued.

\bibliographystyle{elsarticle-num}
\bibliography{Biblist}

\end{document}